\documentclass[12pt,a4paper]{article}
\usepackage{amsmath,amstext,amssymb,amscd, color}
\usepackage{graphicx}

\usepackage{float}
\usepackage{cite}

\usepackage[all]{xy}

\usepackage{tikz}

\usepackage[english]{babel}

\newcommand{\Xcomment}[1]{}

\newtheorem{theorem}{Theorem}[section]
\newtheorem{lemma}[theorem]{Lemma}
\newtheorem{corollary}[theorem]{Corollary}

\newtheorem{prop}[theorem]{Proposition}
\newtheorem{remark}[theorem]{Remark}

\newtheorem{dfn}[theorem]{Definition}

\makeatletter \@addtoreset{equation}{section} \makeatother

\newenvironment{proof}{\noindent{\bf Proof}~}%
{\hfill$\qed$\medskip}

\def\qed{ \ \vrule width.1cm height.3cm depth0cm}

 \makeatletter
\renewcommand{\section}{\@startsection{section}{1}{0pt}%
{-3.5ex plus -1ex minus -.2ex}{2.3ex plus .2ex}%
{\normalfont\Large}}
 \makeatother

 \makeatletter
\renewcommand{\subsection}{\@startsection{subsection}{2}{0pt}%
{-3.0ex plus -1ex minus -.2ex}{-1.5ex plus .2ex}%
{\normalfont\large\bf}}
 \makeatother

 \Xcomment{
 \makeatletter
\renewcommand{\subsection}{\@startsection{subsection}{2}{0pt}%
{-3.0ex plus -1ex minus -.2ex}{1.5ex plus .2ex}%
{\normalfont\large\bf}}
  }
 \makeatother

\def\BN{{\mathbb N}}
\def\Zset{{\mathbb Z}}
\def\BZ{{\mathbb{Z}}}

\def\Bscr{{\cal B}}
\def\CB{{\cal B}}
\def\BBscr{{\cal{BB}}}
\def\CBB{{\cal{BB}}}

\def\CM{{\cal M}}

\def\CO{{\cal O}}

\def\CW{{\cal W}}

\def\tilde{\widetilde}
\def\hat{\widehat}
\def\bar{\overline}

\makeatletter \@addtoreset{equation}{section} \makeatother

 \makeatletter
\renewcommand{\section}{\@startsection{section}{1}{0pt}%
{-3.5ex plus -1ex minus -.2ex}{2.3ex plus .2ex}%
{\normalfont\Large}}
 \makeatother

 \makeatletter
\renewcommand{\subsection}{\@startsection{subsection}{2}{0pt}%
{-3.0ex plus -1ex minus -.2ex}{-1.5ex plus .2ex}%
{\normalfont\normalsize\bf}}
 \makeatother

  \Xcomment{
 \makeatletter
\renewcommand{\subsection}{\@startsection{subsection}{2}{0pt}%
{-3.0ex plus -1ex minus -.2ex}{1.5ex plus .2ex}%
{\normalfont\normalsize\bf}}
 \makeatother
 }

 \makeatletter
\renewcommand{\subsubsection}{\@startsection{subsubsection}{2}{0pt}%
{-2.0ex plus -1ex minus -.2ex}{-2.0ex plus .2ex}%
{\normalfont\normalsize\underline}}
 \makeatother

\def\BR{{\mathbb R}}
\def\Zset{{\mathbb Z}}
\def\BZ{{\mathbb Z}}

\def\Bscr{{\cal B}}

\def\CP{{\cal P}}

\def\tilde{\widetilde}
\def\hat{\widehat}
\def\bar{\overline}

\def\lra{\longrightarrow}
\def\iso{\overset\sim\lra}
\def\isom{\overset\sim{=}}

\begin{document}

\title{Bruhat operads }

\author{Gleb Koshevoy and Vadim Schechtman, with Appendix by Daria Poliakova}

\date{\today}

\maketitle

{\em \qquad \qquad  \qquad  \qquad  \qquad To the memory of Yuri Ivanovich Manin}

\

\

{\em \qquad \qquad  \qquad  \qquad  \qquad M. Jourdan: Par ma foi, il y a plus de quarante ans 

\qquad \qquad \qquad \qquad \qquad  que je dis de la prose sans que j'en susse rien

\ 

\qquad \qquad \qquad \qquad \qquad \qquad \qquad  Bourgeois gentilhomme} 

\

\begin{abstract}

\

We describe some planar operads built from the higher Bruhat orders and show that they admit 
a multiplication.

\end{abstract}

\section{Introduction}

\subsection{} The present paper consists of three parts. The aim of the first part
\newline (Sections \ref{plan-op} - \ref{op-higher-bruhat}) is to introduce some planar operads 
builtbuild from higher Bruhat orders, cf. \cite{MSFA}. We call them {\em Bruhat operads}.  

In the second part (Sections \ref{mult-cos-shifted} - \ref{ursa-major}) we show that the Bruhat operads are {\em operads 
with multiplication} in the sense of \cite{MS}. 

The third part is Appendix by Daria Poliakova with the proof that the Bruhat operad is compatible with Bruhat orderings.

\

\subsection{} Recall that a {\em planar} (sometimes called {\em nonsymmetric}) operad of sets $\mathcal O$ consists of a family of sets 
$\{\mathcal O(n)\}_{n\in \mathbb{Z}, n\geq 0}$, together with composition maps
\[
\gamma:\ \mathcal{O}(n)\times \mathcal{O}(n_1)\times\ldots\times \mathcal{O}(n_m)\to 
\mathcal{O}(\sum_i n_i) 
\]
and a unit $1\in \mathcal{O}(1)$
satisfying the associativity and the unit axiom as for a usual operad, except that we do not require an action of symmetric groups $S(n)$ on $\mathcal{O}(n)$, cf. 
\cite{M}, Def. 3.12 and  Def. \ref{def-planar} below.

\

\subsection{}
Let $d, n$ be integers, $1\leq d\leq n$. In \cite{MSFA} (cf. also \cite{MS}) certain posets $B(n,d)$ 
called {\it higher Bruhat orders} have been introduced,  
their definition is recalled below in \ref{def-hbo}. 

\

One can describe them roughly as follows. $B(n,1) = S(n)$ with  
the classical weak Bruhat order. The poset $S(n)$ admits unique maximal and minimal elements, namely, the minimal element is the unit permutation $e$, whereas the maximal one 
is the permutation $\sigma_{max}$ of maximal length\footnote{where length means the usual length $\ell$ in a Coxeter group, so $\ell(\sigma_{max}) = n/(n-1)/2$.}. By definition 
$B(n,2)$ is a certain quotient of the set of maximal (i.e. noncompressible\footnote{which means a collection $e < x_1 < \ldots < x_n < \sigma_{max}$ which cannot be included into a longer collection}) chains in $S(n)$ connecting $e$ with $\sigma_{max}$. 

\

Afterwards one proceeds by induction on $d$. Each $B(n,d)$ is a poset admitting unique minimal and maximal elements, say $b_{min}^d$ and $b_{max}^d$, and $B(n,d+1)$ is a certain quotient of the set of maximal chains in $B(n,d)$ connecting $b_{min}^d$ with $b_{max}^d$.   

\

A remarkable paper by Ziegler \cite{Z} contains a number of deep constructions and results 
on higher Bruhat orders. They are essential in our study.

\subsection{} Fix a positive integer $d$. Denote $\mathcal{O}(n) := B(nd, d)$. The first main result of the present paper, see Thm \ref{main-thm} below, says that the collection $\mathcal{O} = \{\mathcal{O}(n)\}$ 
admits a structure of a planar operad, to be called the {\em small Bruhat operad}. 

\

Our definition of multiplications $\gamma$ is 
based on the following main construction.

\

Given two elements $b\in B(n,d)$ and $b'\in B(m,d)$ and an integer $j$, $0\leq j\leq 
n+m-d$, we define a new element $b\circ_j b'\in B(n+m-d)$, see \ref{ins-higher-b}.

\ 

The definition of the operations 
$\circ_j$ is nontrivial and uses an important theorem by Ziegler, see Thm. \ref{Ziegler} below. 

\

\subsection{} We may include into the game elements of $B(m,d)$ with $d$ not dividing $m$ as well. Namely, we define a planar operad $\mathcal{O}' = \{\mathcal{O}'(n)\}$ which contains $\mathcal{O}$ as a suboperad and called the {\em big Bruhat operad}. The sets 
$\mathcal{O}'(n)$ are certain infinite unions of sets $B(m,d)$ with fixed $d$ and arbitrary $m \geq d$.

\

This is the second main result of the present paper, see \ref{main-thm-b}.

\subsection{} The contents of the paper is as follows. In Section 2 we recall the definition of planar operads. In Section 3 we discuss certain "Master  planar operads" 
which serve as a combinatorial base for constructing the big Bruhat operad. In Section 4 we introduce the main operation of insertion. In Section 
5 we discuss some operads related to symmetric groups which correspond to the case $d = 1$ 
of our considerations. Section 6 contains the definitions of Bruhat operads.

Sections 7 - 10 are devoted to multiplicative structures. 

\subsection{•} The results if the present paper  have been announced in \cite{KS1}, \cite{KS2}. 

\subsection{Acknowledgement} We are grateful to Jim Stasheff for useful remarks, to Denis Bashkirov for an interesting correspondence, to Vladimir Dotsenko, Alexander Voronov, and Ralph Kaufmann for interesting discussions, and especially to Dasha Polyakova for 
interesting discussions and corrections.

\subsection{Notation} 


If $X$ is a set, $\CP(X)$ will denote the set of all its subsets.  

For $n\in \mathbb{Z}, n\geq 1$, $[n] := \{1, \ldots , n\}$.
\newline For $a\in \mathbb{R}, 
a + [n] := \{ a+1,\ldots a+n\}$.

\

By ${[n]}\choose d$ we denote the set of $d$ element subsets of $[n]$.

\

\

\centerline{Part I. OPERADS}

\section{Planar operads}\label{plan-op}

\subsection{} Let $(K, \otimes)$ be a tensor category with unit object $e$. 
 Below we will be interested mostly in the category of sets with the monoidal structure given by the cartesian product. 

\begin{dfn}\label{def-planar}
A {\em planar operad} in $K$ consists of objects $\CO(n)$ of $K$ indexed by the natural numbers 
$n\in \Zset, n \geq 1$, with  operations of composition

\begin{equation}\label{op1}
\gamma_{k;n_1, \ldots, n_k}\,\colon\, \CO(k)\otimes \CO (n_1)\otimes \cdots \otimes \CO(n_k)\to \CO(n_1+ \cdots +n_k)
\end{equation}
and a unit map $I: e\lra \CO(1)$.

The composition should be associative which means that  

\begin{equation}\label{com1}
\gamma_{k; n_1, \ldots,  n_k }
\prod_1^k
\gamma_{n_t;l_{m_1}^t, \ldots, l_{m_{n_t}}^t}=
\gamma_{k; \sum_{j=1}^{n_1 }l^1_{m_j}, \ldots,  \sum_{j=1}^{n_k }l^k_{m_j}},
\end{equation}
 
the composition 
\[
\CO(k) = \CO(k)\otimes e\otimes \cdots \otimes e \lra 
\CO(k)\otimes \CO (1)\otimes \cdots \otimes \CO(1)\overset\gamma\lra \CO(k)
\]
should be equal to $\text{Id}_{\CO(k)}$, and the composition 
\[
\CO(n) = e\otimes \CO(n) \lra \CO(1)\otimes\CO(n)\overset\gamma\lra \CO(n) 
\]
should be equal to $\text{Id}_{\CO(n)}$.

\end{dfn}

\subsection{Partial compositions, or insertions.}\label{part-comp} To define compositions (\ref{op1}) it is sufficient to define their 
particular cases, "insertions", or "blowups".

\begin{equation}\label{compos-ins}
\circ_i : \CO(m)\otimes \CO(n) \to \CO(m+n-1)
\end{equation}
where
\[
a\circ_i b = \gamma_{m;1,\ldots,1,n,1,\ldots, 1}(a;I,\ldots, I,b,I,\ldots, I),
\]
$1\leq i\leq m$, 
satisfying two identities:

\

for all $a\in \CO(l), b\in \CO(m), c\in \CO(n)$ we should have
\[
(a\circ_i b)\circ_{i-1+j} c = a\circ_i (b\circ_j c),
\] 
$1\leq i\leq l, 1\leq j\leq m$ 
(associativity), and
\[
(a\circ_i b)\circ_{k-1+m}c = (a\circ_k b)\circ_i c,
\]
$1\leq i < k\leq l$ (commutativity), 
see \cite{LV}, 5.9.4.

\ 

This criterion will be crucial for us.

 \section{Master planar operads}
 
 \subsection{Operad $M$.}
 
 Consider the following collection of sets $M(n)$, $n\ge 1$.
 
 The set $M(n)$ is consists of collections of integer numbers  $\mathbf k:=(k_0, k_1, \ldots, k_n)$, $k_i\ge 0$.
 
 \
 
 We can imagine elements of $M(n)$ as collections of black and white different points on the real line such that 
 first (from left to right) we see $k_0$ white points, then a black point, then $k_1$ white points, 
 then the second black point, then $k_2$ white points, and so on until the $n$-th black point and $k_n$ white points.
 
 \   
 
Let us define a family of maps   
\begin{equation}\label{master1}
\gamma_{n; l_1, \cdots, l_n}: \ M(n)\times M(l_1)\times \cdots\times M(l_n)\to M(\sum_{j=1}^n l_j)
\end{equation}
by the following  rule.

 For  $\mathbf k:=(k_0, k_1, \ldots, k_n)\in M(n)$,  $\mathbf m_1:=(m_0^1, m_1^1, \ldots, m_{l_1}^1)\in M(l_1),\cdots $, 
  $\mathbf m_n:=(m_0^n, m_1^n, \ldots, m_{l_n}^n) \in M(l_n)$, we define 
 \[
 \gamma_{n; l_1, \cdots, l_n}(\mathbf k, \mathbf m_1, \cdots, \mathbf m_n) =  
 \] 
\[
 (k_0+m_0^1, m_1^1, \cdots, m^1_{l_1-1},  k_1+m_{l_1}^1+m_0^2,  m_1^2, \cdots, m_{l_2-1}^2, k_2+m_{l_2}^2+m_0^3, \cdots,\]
  \begin{equation}\label{comp1}
  k_{n-1}+m_{l_{n-1}}^{n-1}+m_0^n,
 m_1^n, \cdots, m_{l_n-1}^n,  k_{n}+m_{l_{n}}^n).
 \end{equation}
 
 Geometrically the composition $\gamma$ means that we insert the whole collection $M(1)$ instead of the first black point in 
 $\mathbf k$, $M(2)$ instead of the second black point in 
 $\mathbf{k}$, $\ldots$, $M(n)$ instead of the $n$-th black point in $\mathbf{k}$, see the picture below.  
 
 The unit element is the tuple $I=(0,0)\in M(1)$.
 
 It is easy to check that the map (\ref{master1}) is associative and $I$ is the unit. 
 Thus  we have

 \begin{lemma}\label{master0}
 The collection  $M = \{M(n), n\ge 1\}$, with the  family of maps (\ref{master1}) is a planar operad of sets.
 \end{lemma}
 $\square$

Note that the sets $M(n)$ are infinite.
 
  \subsection{Operads of molecules $M^d$.}\label{operad-molecules}

For a given positive integer $d$,  and  
  $\mathbf k=(k_0, \cdots, k_n)\in M(n)$, consider the  collection
  \[
  \mathbf k^d := (k_0, d, k_1, d, k_2, d, \cdots , k_{n-1}, d, k_n).
  \]
  The sum of these numbers will be denoted by $|\mathbf{k}^d|$; it is equal to $dn+\sum_j\ k_j$. 
  
  \
  
  A sequence $\mathbf{k^d}$ will be called {\em a type}.

  \
  
  Let us consider  the   set $M(\mathbf{k}^d) = [dn +\sum_j k_j]$. Its elements will be called {\em particles}.   
  
  \ 
  
  Subsets of $M(\mathbf{k}^d)$ having the form  of intervals 
  $$
  k_j + dj + [d] = \{k_j+dj+1, \cdots, k_j+d(j+1)\}
  $$  
  we will call {\em teams}, or {\em nuclei}, $j=0, \cdots, n-1$. Elements of nuclei will be called 
  {\em protons}; so each nucleus contains $d$ protons. 
  
  \ 
  
  Particles of $M(\mathbf{k}^d)$ which do  not belong to any nucleus will be called {\em idles}, or {\em electrons}. 
  
  \
   
  The whole set $M(\mathbf{k}^d)$ equipped with the collection of 
  $n$ nuclei inside it will be called a {\em molecule of type $\mathbf{k}^d$}. Thus $M(\mathbf{k}^d)$ contains $|\mathbf{k}^d|$ particles, among them $nd$ protons and 
  $\sum_j k_j$ electrons.
  
  \ 
  
  Obviously the type defines the molecule uniquely.  
  

\ 

Denote by $M^d(n)$ the set of all molecules with $n$ nuclei and arbitrary number of electrons. We can identify $M^1(n)$ with $M(n)$, with nuclei (resp. electrons) corresponding to black (resp. white) points.

\     
   
Define a collection of maps

\begin{equation}\label{masterd}
\gamma_{n; l_1, \ldots, l_n}:\ M^d(n)\times M^d(l_1)\times \cdots\times M^d(l_n)\to M^d(\sum_{j=1}^n l_j)
\end{equation}

 For molecules $M(\mathbf{k})\in M^d(n)$,  $\mathbf k =(k_0, d,  k_1,d, \ldots,d,  k_n)$, 
 $M(\mathbf m_1)\in M^d(l_1)$,   
 $\mathbf m_1 =(m_0^1,d,  m_1^1, d, \ldots, d, m_{l_1}^1),\ldots $, $M(\mathbf m_n)\in M^d(l_n)$,  
  $\mathbf m_n = (m_0^n, d, m_1^n, d, \ldots, d,  m_{l_n}^n)$, we set
 
 \[
 \gamma_{n; l_1, \ldots, l_n}(
 M(\mathbf k), M(\mathbf m_1), \ldots,  M(\mathbf m_n)) = 
 \] 
\[
 M(k_0+m_0^1, d, m_1^1, d, \ldots, m^1_{l_1-1}, d, k_1+m_{l_1}^1+m_0^2,  d, m_1^2, \ldots, m_{l_2-1}^2,  d, k_2+m_{l_2}^2+m_0^3, d, \ldots,
 \]
  \begin{equation}\label{compd1}
 d, k_{n-1}+m_{l_{n-1}}^{n-1}+m_0^n, d, 
 m_1^n, \ldots, d ,m_{l_n-1}^n,  d, k_{n}+m_{l_{n}}^n)).
 \end{equation}
 
\begin{figure}[H]
\centering
\includegraphics[width=14cm, height=14cm]
{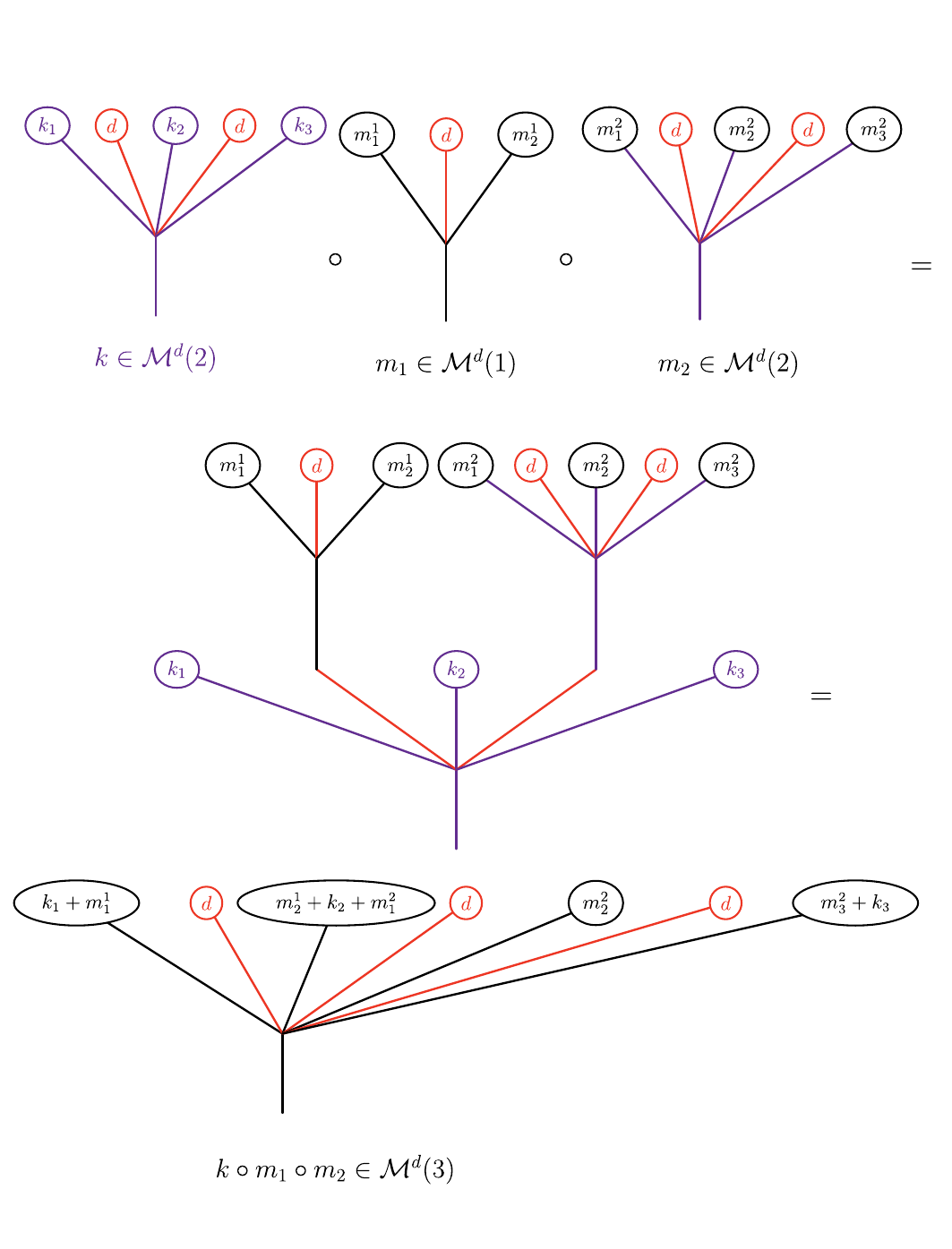}
\caption{Composition in the Master operad}
\end{figure}
  
  \
  
  This operation means that instead of the $j$-th nucleus in $M(\mathbf{k})$ we insert the molecule 
  $M(\mathbf{m}_j)$, $1\leq j\leq n$. 
  
  \

\begin{lemma}\label{oper-md}
The collection   $M^d = \{M^d(n), n\ge 1\}$, with the family of maps (\ref{compd1})  is  a planar operad  
isomorphic to $M^1$. 
\end{lemma}
  
   {\em Proof} Obvious.  $\Box$
   
   \
 
 Note that the sets $M^d(n)$ are infinite.

\section{Symmetric groups}

\subsection{A symmetric operad $S$.}
It is  wellknown that the  symmetric groups form an operad.

Namely, let $S(n)$ be the group of permutations of the set $[n]=\{1, \cdots, n\}$.

The family of maps
\begin{equation}\label{symmop}
\gamma_{k,n_1, \cdots, n_k}: S(k)\times S(n_1)\times \ldots \times S(n_k)\to S(n_1+\ldots +n_k)
\end{equation}
defined by sending

\[
(\tau, \sigma_1, \ldots \sigma_k), \quad \tau\in S(k), \, \sigma_i\in S(n_i), \, i=1, \cdots , k,
\]
to the  permutation of the set $[n_1+\cdots +n_k]$
obtained by applying $\sigma_i$ to the $i$th interval $n_1+\cdots +n_{i-1}+1, \ldots , n_1+\cdots +n_{i-1}+n_i$, and then permuting these intervals following $\tau$.

 The collection $S = \{ S(n), n\ge 1\}$, with the above family of maps $\gamma$ is a planar operad, and even more  it is a symmetric operad.

\subsection{A planar operad $F$.}
 
 We introduce another planar operad $F$ on symmetric groups which is related to  the master operad $M^1$.
 
 We define $F(n)$ as a subset of the set  of  permutations  of $M^1(n)$. 
 
 Specifically,  we consider  collections 
 \[
\mathbf{k} = (k_0, 1, k_1, 1, k_2, \cdots, k_{n-1}, 1, k_n)\in M^1(n)
 \] 
 with teams of size $1$.

 For any set $[n+\sum_{i=0}^n k_i]$ we consider  permutations  $\sigma\in S(n+\sum_ik_i)$ which  permute only nuclei and are identical on the electrons. 
 
 The set $F(n)$ consists of all such permutations for all  $\mathbf k\in  M^1(n)$.

 Compositions  
 
 \[
 \gamma: F(k)\times F(n_1)\times \ldots \times F(n_k)\to F(n_1+\ldots +n_k)\]
 are defined as follows.
 
 \
 
 For $\tau\in F(k)$, $\sigma_i\in F(n_i)$, $i=1, \ldots, k$,
 \[
 \gamma(\tau; \sigma_1, \ldots , \sigma_k)
 \]
 permutes 
 \[
 \gamma(\mathbf k; \mathbf m_1, \ldots, \mathbf m_k), \quad \mathbf k\in M^1(k), \, \mathbf m_i\in M^1(n_i), i=1, \ldots, k,
 \]
 by permuting the blocks $\mathbf m_i$  by $\sigma_i$ and then permuting these blocks by $\tau$.
 
 We have 
 \begin{prop}
 The collection of $F(n)$, $n=1, \cdots$ with the above defined family of maps is a planar operad.
 \end{prop}

 {\em Proof} follows from  Lemma \ref{oper-md} for $d = 1$.
 \hfill $\Box$

\begin{remark}
Restricting the operad $F(n)$ to the subset of collections without idles gives us the operad $S$ on symmetric groups $S(n)$, $n=1, \ldots$.
\end{remark}

\

  \subsection{Weak Bruhat order}
  
On a symmetric group $S(n)$, there is  a weak Bruhat order $\prec$ which may be defined as follows: 
$u\prec v$ iff the sets of inversions are nested, $Inv(u)\subset Inv(v)$.

\

For example the weak Bruhat order for $S(3)$ looks as follows: 
\newline $e < (12) < (123) < (13)$ and 
$e < (23) < (132) < (13)$.

\ 

 Both operads $S$ and $F$ are compatible with weak Bruhat orders. Namely
 
\begin{theorem}\label{order-1}
 Let $\tau, \tau'\in S(k)$ (resp. $\tau, \tau'\in F(k)$), $\sigma_i, \sigma'_i\in S(n_i)$ 
 (resp. $\sigma_i, \sigma'_i\in F(n_i)$), $i=1, \ldots k$. If $\tau\prec\tau'$ and  
$\sigma_i\prec\sigma'_i$, $i=1, \ldots, k$, then
  \[
 \gamma(\tau; \sigma_1,\ldots , \sigma_k)\prec 
 \gamma(\tau'; \sigma_1',\ldots,  \sigma_k').
 \]
\end{theorem}

We prove this theorem below in Section \ref{invsec}, see \ref{proof-order}.

\


\section{Higher Bruhat orders and insertions}\label{invsec}

\subsection{Higher Bruhat orders.}\label{def-hbo}

Let us recall some definitions from \cite{MSFA}. 

\

Let $d \leq n$ be two positive integers. Let $Lin(n,d)$ denote the set of linear orders on the discrete Grassmannian ${[n]}\choose d$. 

\ 

For a subset $\mathbf{i} = \{i_1, \ldots, i_{d+1}\}\in {{[n]}\choose d+1}$, 
$1\leq i_1<\cdots < i_{d+1}\leq n$, and $1\leq k\leq d+1$ let us denote $\hat{\mathbf{i}}_k$ 
the sequence $\mathbf{i}$ with $i_k$ omitted. 

\

Let us call its {\it packet}, and denote by $P(\mathbf{i})$, the sequence
\[
P(\mathbf{i}) = (\hat{\mathbf{i}}_{d+1}, \hat{\mathbf{i}}_d, \ldots, \hat{\mathbf{i}}_1).
\]      

\

Let us call an element $\prec\in L(n,d)$ {\it admissible} if 
its restriction to each packet $P(\mathbf{i})$ is 
either {\em lexicographical}, i.e.
\[
\hat{\mathbf{i}}_{d+1}\prec \hat{\mathbf{i}}_d \prec \ldots \prec \hat{\mathbf{i}}_1
\]
or {\em antilexicographical}, i.e.
\[
\hat{\mathbf{i}}_{d+1}\succ \hat{\mathbf{i}}_d \succ \ldots \succ \hat{\mathbf{i}}_1.
\]
Let $A(n,d)\subset L(n,d)$ denote the subset of admissible linear orders. 

\

The set $B(n,d)$ of Bruhat 
orders for ${[n]}\choose d$  is defined as the quotient $A(n,d)/\sim$ where $\sim$ is an equivalence relation described in 
\cite{MSFA}, no 2.

\

Each $B(n,d)$ is a poset, to be denoted $\mathcal{HB}^d(n)$. One of main results of  \cite{MSFA} states that any order from $\mathcal{HB}^d(n)$ can be included into a maximal chain joining the minimal order with the maximal one. 

\subsection{Example} $A(n,1) = B(n,1) = S(n)$ with the classical weak Bruhat order.

\

\subsection{Remark} Alternatively $B(n,d)$ may be defined as the set of minimal  
{\em partial} orders on ${[n]}\choose d$ whose restriction to each packet is either 
lexicographical or antilexocographical.

\
{
Recall that partial orders, anti-symmetric transitive binary relations, are partially ordered by the inclusion of ordered pairs of elements. If we consider 
the intersection of all partial orders that have the same restriction to every packet, we get a  minimal  
{\em partial} orders on ${[n]}\choose d$ whose restriction to each packet is either 
lexicographical or antilexicographical.
}

\

\subsection{Inversion sets} 
Let $a\in A(n,d)$. The 
collection of all $\mathbf{i}\in$ ${[n]}\choose {d+1}$ such that $P(\mathbf{i})$ is ordered anilexicographically forms the {\em set of inversions } $Inv(a)\subset$ 
${[n]}\choose {d+1}$. If $a\sim a'$ then $Inv(a) = Inv(a')$, so 
for $\prec\in B(n,d)$ the set 
  $Inv(\prec)\subset$ ${[n]}\choose {d+1}$ is well defined.
  
This way we get a map

\

$Inv:\ A(n,d) \lra \CP$ ${[n]}\choose {d+1}$

\

The set of inversions determines the order uniquely. So $B(n,d)$ may be defined as the image 
of $Inv$. 

\

\subsection{Ziegler criterion} 
Ziegler \cite{Z}, Th. 4.1(B), characterised the subsets of ${[n]}\choose {d+1}$ which appear as sets of inversions of  Bruhat orders.

\

\begin{theorem}\label{Ziegler} {\em (Ziegler)}
A subset $I\subset {{[n]}\choose {d+1}}$ is the set of inversions $Inv(\prec)$ of a  Bruhat order $\prec$ for ${[n]}\choose d$ if and only if for any $\mathbf{i} \in {{[n]}\choose {d+2}}$
the intersection of $I$ and the packet $P(\mathbf{i})$  is either a beginning or an ending interval of $P(\mathbf{i})$.
\end{theorem}

\

The last condition means that there exists $k$, $0\leq k\leq d+2$ such that 
$P(\mathbf{i})\cap I$ is equal to 
$(\hat{\mathbf{i}}_k, \hat{\mathbf{i}}_{k-1},\ldots, \hat{\mathbf{i}}_1)$
or to $(\hat{\mathbf{i}}_{d+2}, \hat{\mathbf{i}}_{d+1},\ldots, \hat{\mathbf{i}}_k)$. Here $k=0$ (resp. 
$k = d+2$) means that $P(\mathbf{i})\cap I = \emptyset$ (resp. $P(\mathbf{i})\cap I = P(\mathbf{i})$).

\
 
We will say that {\em $\mathbf{i}$ satisfies the Ziegler criterion}, or simply that {\em $\mathbf{i}$ satisfies} (Z), for $I$.

\

Let us denote by $Z(n,d+1)$ the set of subsets $I\subset {{[n]}\choose {d+1}}$ such that for any 
$\mathbf{i} \in {{[n]}\choose {d+2}}$, $\mathbf{i}$ satisfies (Z) for $I$. 
The above statements mean that $Inv$ establishes a bijection
\[
Inv:\ B(n,d) \overset\sim\to Z(n,d+1).
\]   

\

Note that a maximal chain in the poset $\mathcal{HB}^d(n)$  has the following chain of inversion sets. It starts with the empty set, and each consequent one is obtained by adding one more inversion. The length of such a chain is $[n]\choose d+1$, the cardinality of the inversion set for the maximal element of $\mathcal{HB}^d(n)$.

\subsection{Topological definition: wire spaces $\CW(n,d)$.}

\subsubsection{\CW(n,1)} We define the wire space 
$\CW(n,1)$ as the topological space whose elements are collections 
of distinct real points $x_1 < x_2 <\ldots < x_n$ on the real line together with a bijection 
$\sigma:\ [n] \iso [n]$. We assign to $x_i$ the number $\sigma(i)$. 

So $\CW(n,1)$ is a subspace of $(\BR^d)^{n!}$, and we consider it as equipped with the induced topology. We have an obvious bijection
\[
\pi_0(\CW(n,1))\isom S(n) = B(n,1), 
\]
and all connected components of $\CW(n,1)$ are contractible. 

\subsubsection{\CW(n,2)}. Elements of this space are "arrangements of $n$ pseudolines in 
$\BR^2$".

More precisely, an element of $\CW(n,2)$ is the following collection of data. 

\

(a) "Walls": two distinct vertical lines $w_l = \{ x_l + iy, y\in \BR\}, 
w_r = \{ x_r + iy, y\in \BR\}, y_l < y_r$   in $\BR^2$, each of them equipped with a collection of $n$ distinct points:
\
\[
a_i = x_l + iy_j, y_n > \ldots > y_1,\ 
\]
\[
b_i = x_r + iz_j, z_1 > \ldots > z_n.
\]

\

(b) "Wires": $n$ wires, the $j$-th wire $w_j$ being the graph of a smooth function connecting 
$a_j$ with $b_j$, 
\[
f_j: \ [x_l, x_r] \lra \BR, f(x_l) = y_j,\ f(x_r) = z_j.
\]
We require that for all $j < k$ the wire $w_j$ meets $w_k$ exactly once.

\subsection{Theorem.} {\em We have a canonical bijection 
\[
\pi_0(\CW(n,2))\isom B(n,2)
\]
and all connected components of $\CW(n,2)$ are contractible.} $\Box$

\
{
Proceeding in the same manner one can give a similar description of the higher 
Bruhat sets $B(n,d)$  with pseudohyperplanes and get \cite{VK,Z}
the bijection \[
\pi_0(\CW(n,d))\isom B(n,d),
\]
where $\CW(n,d)$ denote the collections of $n$ pseudohyperplanes in $\mathbb{R}^d$, such that any $d$ of them intersect once in a single point.    
 }

 \subsection{Main construction: insertions of higher Bruhat orders}\label{ins-higher-b}
 
 To define an operad structure on higher Bruhat orders we use the following construction.
 
 \
 
 Let $d, n, m$ be positive integers, $n > d, m > d$, and $j$ an integer such that 
 $$
 0\leq j\leq n+m-d.
 $$ 
 Given two Bruhat orders $\prec$ for ${[n]}\choose {d}$ and $\prec'$ for ${[m]}\choose {d}$ we will define a new Bruhat order $\prec\circ_j\prec'$ for ${[n+m-d]}\choose {d}$, to be called 
 {\em the $j$-th insertion of a 
   $\prec'$  into  $\prec$}. 
 
 
\ 

We define such an insertion by specifying  its  set of inversions.
 
 \
 
  Namely we define a subset $I(\prec\circ_j\prec')\subset$ ${[n+m-d]}\choose {d+1}$ and afterwards prove that it satisfies (Z), whence it uniquely defines an order $\prec\circ_j\prec'$. 
 
 \
 
  Let $Inv(\prec)$ and $Inv(\prec')$ be the sets of inversions of $\prec$ and $\prec'$, respectively.
 
 
 Consider a decomposition 
 \[
 [n+m-d]=[ j]\cup ( j+[m]) \cup (j+m+[n-(j+d)]).
 \]
 We denote $A:=[j]$, $B:=j+[m]$, $D = j + [d]\subset B$, and $C:=j+m+[n-(j+d)]$.
 
 \

(i) First,   for a $d+1$ element subset  $J=\{j_1, \cdots, j_{d+1}\}\subset B$, we say that 
\begin{equation}\label{inv1}
J\in I (\prec\circ_j\prec' )\mbox{ iff  }
  J-j=\{j_1-j, \cdots, j_{d+1}-j\}\in Inv(\prec').
  \end{equation}
  
  
  \

(ii) Second,  for a  $d+1$ element  subset 
\[
I=\{i_1, \cdots, i_{d+1}\}\subset A\cup D\cup  C
\]

\

we say that 
  \begin{equation}\label{inv2}
  I\in I(\prec\circ_j\prec') \mbox{ iff } \phi^{-1}(I)\in Inv(\prec),
 \end{equation}
  where
   \[
   \phi: [n] \to A\cup D \cup C
   \]
   is the unique monotone bijection.
   
   \
   
 (iii)   Third, let $L=\{l_1, \cdots, l_{d+1}\} \subset [n+m-d]$ be none of types (i) or (ii). 
  Denote 
  \[
  L^+=L\cap (A\cup D\cup C)
  \]
  and  
  \[ 
  L^-=L\cap (B\setminus D) = L\setminus L^+.
  \] 
 
 \
 
 
 Thus $L = L^-\coprod L^+$.
 
 \
 
 Let   
 \[f:L^-\to D\]
 be the unique monotone map sending $L^-$ to the rightmost elements of $ D\setminus L^+$. 

\

Denote by 
 \[
 \overline L=L^+\cup f(L^-)\subset A\cup D\cup C
 \]
 (the union is disjoint).
 
 \
 
 {\bf Example.} If $L\cap (j+[m])=\{a, b\}$,  $a\le j+d$, $b>j+d$, then 
 \[
 \overline L=\begin{cases}
  L\setminus b\cup \{j+d\}& if \quad a\neq j+d\cr
 L\setminus b\cup \{j+d-1\} & if\quad  a=j+d.
 \end{cases}
 \] 
 
 $\Box$  
 
 \

 Then we say that 
  \begin{equation}\label{inv3}
  L\in I (\prec\circ_j\prec' ) \mbox{  iff }  \phi^{-1}(\overline L) \in Inv(\prec).
  \end{equation}
  
  This three step construction defines a set  $I(\prec\circ_j\prec')$.
  
  \

By definition
\begin{equation}\label{coprod-I}
I (\prec\circ_j\prec') = I_1\coprod I_2\coprod I_3
\end{equation}
where $I_1$ (resp. $I_2$, $I_3$) are the subsets of elements satisfying (\ref{inv1}) (resp. (\ref{inv2}), 
(\ref{inv3})). 
 
   \begin{prop}\label{inv-prop}
 The set $I(\prec\circ_j\prec')$ defined above satisfies (Z), whence $I(\prec\circ_j\prec') = Inv(\prec\circ_j\prec')$ for the uniquely defined $\prec\circ_j\prec'\in B(n+m-d,d)$.
 \end{prop}
 
 This will be proven below, see \ref{ziegler-d}.

\subsection{Insertion of wiring diagrams}
Let $W$ be a wiring diagram with $m$ wires and $U$ be a diagram with $n$ wires (strands). Let us describe the insertion $W\circ_j U$ following our construction of the insertion 
$B(m,2)\circ_{j} B(n,2)$.

Pick the intersection point $v$ of  wires $j+1$ and $j+2$, and consider a small disc $D_{j+1,j+2}$ with  center at $v$. Let us orient the boundary of this disc clockwise. Then the wire $j+1$ first enters and leave the disc, and and $j+2$ second enters and leave.

We do the insertion   in two steps.

First, we place $W$ inside a disc $D_U$ with the clockwise  oriented boundary such that beginning points of wires are clockwise $1$, $2$, etc and leaving points are clockwise $1$, $2$, etc.
Then we replace the disc $D_{j+1,j+2}$ with the center at $v$ by the disc which contains $U$ we identify the entering and leaving points of wire $j+1$ with that of $1$ and $j+2$ with $n$ (we shift numeration by adding $j$ to indices of $[n]$ and shift by $n+j$, the indices  $\{j+3, \cdots, m\}$ (they are outside of $U$, as well as wires with indices $\{1, \cdots, j\}$).

Second, outside the disc, we extend the wire $1$ of $ U$ as $j+1$ in $W$, and the wire $m$ with $j+2$, and wires $2, \cdots, m-1$ we etend to be parallel to $j+2$  in $W$ and parallel between themselves. 

\subsection{Remark.} The above insertion resembles the little disc operad, but we have to expand wires outside inserted discs while in the little disc operad one step of insertion is enough.

\subsection{}
Let us show that the above insertion gives the same result as \ref{ins-higher-b}.  

First, we recall meaning of the inversion in wiring diagrams. Consider three wires with indices $i<j<k$. Then this triple form an inversion if $i$ first intersect $k$, and no inversion if $i$ first intersect $j$.

Now, in the above insertion $W\circ_j U$,

\begin{itemize}
\item
 triples of wires with indices of $\{j+1, \cdots, j+m\}$ have the same   intersection as in $U$, that is (5.1);
 
\item triples of wires with indices different of $\{j+1, \cdots, j+m\}$ have the same   intersection as in $W$, that is (5.2);

\item a triple of wires with one index of  $\{j+1, \cdots, j+m\}$ has the same inversion as this index is $j+2$  and the corresponding intersection is as in $W$, that is as in (iii);
\item a triples with two indices of  $\{j+1, \cdots, j+m\}$ has the same intersection as we regard these two indices $j+1$ and $j+2$ in $W$. In fact, $j+1$ and $j+2$ are parallel outside the disc $D_{j+1,j+2}$., that is as in (iii).

\end{itemize}

\subsection{Remark.} {In dimension $d\ge 2$ the items (5.1) and (5.2) can be interpreted geometrically in the language pseudohyperplane arrangements for $B(n,d)$.

\


\

Namely, we remove pseudohyperplanes with indices $j+1, \cdots , j+d$ from $W$. In a small $d$-dimensional ball which has a center at the intersection of pseudohyperplanes with indices $j+1, \cdots , j+d$ in $W$, we insert the pseudohyperplane arrangement $U$. Then we extend $U$ outside the ball following the step (iii)  the insertion construction. 
However, we do not have a geometric interpretation that, 
beyond the case $d=2$. See also \cite{Bash}.

} 
 \subsection{Example for $d = 2$.} The picture below gives an example of an insertion.   
 
We start with two elements $w_1\in B(5,2)$ drawn as a collection of $5$ wires,  and  $w_2\in B(4,2)$ drawn as a collection of $4$ wires. 
 
 \
 
 In the picture below the element 
 $w_1\circ_2 w_2\in B(7,2)$ is drawn.
 
 \
 
 
 \
 
 
 \newpage
 
\begin{figure}[H]
\centering
  \includegraphics[width=14cm, height=14cm]
{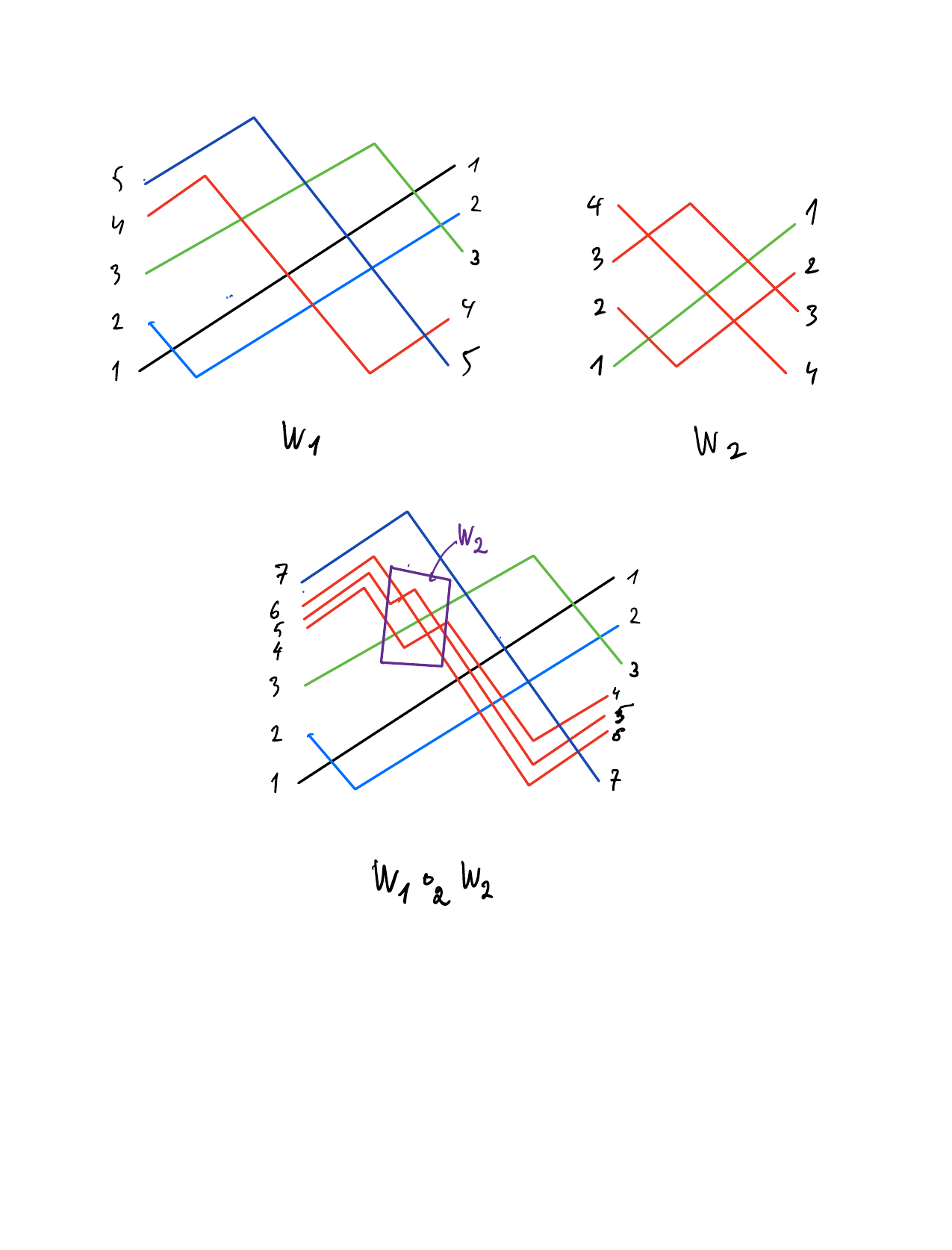}
\caption{An example of an insertion.}
\end{figure}

\

 
 \begin{corollary}\label{compat}
 Let $\prec'$ be dominated by $\prec''$, i.e. \\
 $Inv(\prec')\subset Inv(\prec'')$.
 Then $Inv(\prec\circ_j\prec')\subset Inv(\prec\circ_j\prec'')$.
 \end{corollary}
 
 {\bf Proof}. Let 
 \[
 Inv (\prec\circ_j\prec'') = J_1\coprod J_2\coprod J_3
 \]
 be the decomposition (\ref{coprod-I}) for $Inv (\prec\circ_j\prec'')$. Then $I_1\subset J_1$, whereas 
 $I_2 = J_2$ and $I_3 = J_3$.   
 \hfill $\Box$
 
 \begin{corollary}\label{compat2}
 Let $\prec$ be dominated by $\tilde\prec$, i.e. \\
 $Inv(\prec)\subset Inv(\tilde\prec)$.
 Then $Inv(\prec\circ_j\prec')\subset Inv(\tilde\prec\circ_j\prec')$.
 \end{corollary}
 
 {\bf Proof} similar to the proof of the previous corollary.
 
 \subsection{Proof of Thm \ref{order-1}}\label{proof-order} follows from the above two corollaries since 
 the map $\gamma$ is a composition of insertions. $\Box$
 
 \
{
For $d=1$ or $d=2$, the inclusion of the inversion sets, $Inv(\prec)\subset Inv(\prec')$  implies the 1- and 2- Bruhat ordering $\prec$ and $\prec'$. 
Because of that, we get that our Bruhat operad is compatible with Bruhat orderings also for $d=2$. For $d\ge 3$, this is 
also true, see Appendix.

} 
 \
 
 
 \subsection{Weak Bruhat orders on symmetric groups}
 
 For a symmetric group $S(k)$, a subset $W\subset \Omega:=\{(i,j), i<j, i, j\in [k]\}$ is an inversion set of a permutation of $S(k)$ if and only if for any $i<j<l$ if $(i,j)$ and $(j,l)\in W$, then $(i,l) \in W$,
 and if $(i,l)\in W$, then at least one of $(i,j)$ and $(j,l)\in W$.
 This is nothing but Theorem \ref{Ziegler}.
 
 Recall that the weak Bruhat order  $\tau\prec \tau'$ is defined  by the inclusion $Inv(\tau)\subset Inv(\tau')$.
 
 \
 

A packet is a pair $(i,j)$, $i<j$ and a weak Bruhat  order is a linear order on $[n]$ such that for each pair $(i,j)$ either 
$i\prec j$, or $i\succ j$. This is nothing but an aritrary permutation of $[n]$. 

\

\subsection{Insertion of permutations.}\label{ins-perm}

Let $\sigma\in S(n), \tau\in S(m)$, $1\leq j \leq n$. 
Let us define $\sigma\circ_j\tau\in S_{n+m-1}$ as follows. 

\

We will use two decompositions
\[
[n] = [j-1]\cup \{ j\}\cup (j+[n-j]) = A\cup \{ j\}\cup B
\]
and
\[
[n] = [\sigma(j)-1]\cup \{\sigma(j)\}\cup (\sigma(j)+[n-\sigma(j)]) = 
A'\cup \{\sigma(j)\}\cup B'. 
\] 
Obviously $\sigma$ induces an isomorphism 
\[
\sigma: A\cup B \overset\sim\to A'\cup B'.
\]
Denote for brevity $j':=\sigma(j)$. 

\

Now let us stretch our decompositions to 
\[
[n+m-1] =  A\cup j - 1 + [m]\cup m + B 
\]
and
\[
[n+m-1] = A'\cup j' - 1 + [m]\cup m + B'. 
\]
As before, $\sigma$ induces an isomorphism
\[
\sigma: A\cup m + B \overset\sim\to A'\cup m + B'.
\]
Moreover, $\tau$ induces an isomorphism 
\[
\tau: j - 1 + [m] \overset\sim\to j' - 1 + [m].
\]

\

Now we define $\rho := \sigma\circ_j\tau\in S(n+m-1)$ as follows: 

\

$\rho(x) = \sigma(x)$ if $x\in A\cup m + B$, 

\

$\rho(x) = \tau(x)$ if $x\in j - 1 + [m]$.  

\

\subsection{Proof of Proposition \ref{inv-prop} for $d = 1$.} For  weak Bruhat orders $\prec$ and $\prec'$ on $[n]$ and $[m]$, let $\tau$ and $\tau'$ be the  corresponding permutations. 

\

Due to the above construction, for $j\in [n]$, the set  $Inv (\prec\circ_j\prec')$ consists of 
pairs $(j+l, j+k)$, $l<k\le j+m$, if $(l,k)\in Inv (\tau')$ or pairs $(l, j+s)$, $l<j$, $s\le m$, if $(l,j)\in Inv(\tau)$, or $(j+t, k)$, $t\le m$, $j+t<k$ if $(j, k-m+1)\in Inv(\tau)$.

\

But such a set of inversions is exactly the inversion set of $\tau\circ_j\tau'$. $\Box$

\
 

 \subsection{Proof of Proposition \ref{inv-prop}: arbitrary $d$.}\label{ziegler-d} 
 
 We demonstrate that  the Ziegler criterion holds for the above defined set $I(\prec\circ_j\prec')$.
 
 \
 
 Let $L = \{l_1, \cdots , l_{d+2}\}\subset [n+m-d]$, $l_1 < \ldots < l_{d+2}$. We use the notations of 
 \ref{ins-higher-b}.
 
 \
 
 (i) If $L\subset B$ then the Ziegler condition for the packet $P(L)$ holds true due to (\ref{inv1}).
 
 \
 
 (ii)  If $L \subset A \cup D\cup C$ then the Ziegler condition for the packet $P(L)$ holds true because of (\ref{inv2}).
 
 \
 
 (iii) Now let  $L$ be not of the above two types. We define as before 
 \[
 L^+ = L\cap (A\cup D\cup C),
 \]
 \[
 L^- = L\setminus L^+ = L\cap (B\setminus D).
 \] 
 
Thus $L=L^+\cup L^-$.

\

 The packet  $P(L)$ takes the form 
 \[
 P(L) = \{ L\setminus l_{d+2}, \cdots, L\setminus l_1\} = \{L_{d+2}, \ldots, L_1\}
 \]
 
 
 \
 
(iii-a) Consider first the case $L\cap B=\{l_2, \cdots, l_{d+2}\}$. 

\
 
 Due to (\ref{inv3}) 
 
 \
 
 In  such a case 
 
 \[
 \{\phi^{-1}(L_{d+2}), \ldots, \phi^{-1}(L_{1})\} = 
 \]
 \[
 \{(l_1, j+1, \cdots, j+d), \ldots, (l_1, j+1, \cdots, j+d), L_1\}.
 \]
 Obviously, $P(L)$ satisfies the Ziegler condition whether or not  
 
 $(l_1, j+1, \cdots, j+d)$ belongs to $Inv(\prec)$ and  $L_1$ belongs to $Inv(\prec')$.

 \
 
 (iii-b)
 
 Let 
 \[
 g:\ L^- \to D  
 \]
 be the unique monotone map sending $L^-$ to the rightmost elements of $D\setminus L^+$. 
 Denote by 
 \[
 \overline L=L^+\cup g(L^-)\subset A\cup D\cup C
 \]
 (the union is disjoint).

 For $L\cap C=\{ l_{t+1}, \cdots, l_{d+2}\}$, there holds
 \[\overline{L\setminus l_{j}}= \overline L\setminus l_{j}, \quad j=t+1, \cdots, d+2.
 \]
 
For $L^-=\{ l_{s+1}, \cdots , l_{t}\}$, there holds
  \begin{equation}\label{delet1}
  \overline{L\setminus l_{j}}= \overline L\setminus g(l_{s+1}), \quad j=s+1, \cdots, t.
 \end{equation}
 
 For $L\cap (A\cup D)=\{l_{r+1}, \cdots, l_s\}$,
 we have 
  \[\overline{L\setminus l_{j}}= \begin{cases}
  \overline L\setminus g(l_{s+1}), & if \quad l_j>g(l_{s+1})\cr
   \overline L\setminus l_j, & if  \quad l_j<g(l_{s+1}).
   \end{cases}
 \] 
 
 For $L\cap A=\{l_1, \cdots, l_r\}$, there holds
   \[\overline{L\setminus l_{j}}= \overline L\setminus l_j, \quad j=1, \cdots, r.
 \]
 
Thus, for such $L$,  the sets of the packet $P(L)$  coincide with the sets of the packet $P(\overline L)$
except for $(t-s)$ sets corresponding to deleting from $\overline L$ elements of $g(L^-)$. Namely, sets of $P(L)$ corresponding to deleting elements of $L^-$ from $L$,
are the same, according to (\ref{delet1}), and equal  to   
\[
 \overline L\setminus g(l_{s+1}).\]
 
However, $\overline L\setminus g(l_{s+1})$ is a set of the packet  $P(\overline L)$ and the packet $P(L)$, and  sets of $P(\overline L)$ are the same at the interval  ending interval of $P(L)$ which starts at 
$\overline L\setminus g(l_{s+1})$, and 
 at the beginning interval which ends on the set obtained by deleting $l_{t+1}$. Therefore, the Ziegler condition for $P(L)$ holds true, since the beginning  or ending interval of  inversion sets of $P(\overline L)$  gives also the beginning  or the ending interval of  inversion sets of $P(\overline L)$, and the difference is the repeating set $\overline L\setminus g(l_{s+1})$ in $P(L)$.

\hfill $\Box$

\

More briefly:

\

 (iii-a) If  $d+1$ elements of $L$ are  in $B$ then the packet $P(L)$ has a form 
 \[
 P(L) = (M,\ldots,M, L_1).
 \] 
 Such a packet obviously satisfies (Z). 

\

(iii-b) There are at most $d$ elements of $L$ in $B$. 

Let 
 \[
 g:\ L^- \to D  
 \]
 be the unique monotone map sending $L^-$ to the rightmost elements of $D\setminus L^+$. 
 Denote by 
 \[
 \overline L=L^+\cup g(L^-)\subset A\cup D\cup C
 \]
 (the union is disjoint).

The packet $P(L)$ has a form 

\[
P(L) = 
(\overline L_{d+2}, \cdots, \overline L_{t+1},\overline L_t,\overline L_t, \ldots, \overline L_t, \overline L_{t-s}, \overline L_{t-s-1},\ldots, \overline L_1)
\] 
where $s = Card(L^-\setminus D)$. Obviously such a collection satisfies the Ziegler criterion if the packet  
$P(\overline L)$ satisfies it. $\Box$

\subsection{Associativity and commutativity}

 \begin{lemma}\label{assoc}
 The insertion $\prec\circ_j\prec'$ is associative. That is 
 \[ (\prec\circ_j\prec')\circ_{j+j'}\prec''= \prec\circ_j(\prec'\circ_{j'}\prec'').
 \]
 \end{lemma}
 {\bf Proof.}
 Consider a $d+1$-element subset $L$ of $[n+m+m'-2d]$ obtained by insertion of $[m]$ instead of the team 
 $j+[d]$ into $[n]$ and $[m']$ instead of the team $j+j'+  [d]$ into $[n+m-d]$.
 If $L$ has a non-empty intersection with $A$ or $C$, then (\ref{inv2}) and (\ref{inv3}) gives the same $\overline L$ for $(\prec\circ_j\prec')\circ_{j+j'}\prec''$ and $\prec\circ_j(\prec'\circ_{j'}\prec'')$.
 
 If $L\subset (j+[m])$ then by (\ref{inv1}), $\overline L$ is defined by (\ref{inv1}-\ref{inv3}) for $\prec'\circ_{j'}\prec''$, that is transformed following (\ref{inv1}) or (\ref{inv3}). Both transformation gives the same change of $L$ for 
$(\prec\circ_j\prec')\circ_{j+j'}\prec''  $, since $L$ does not change in $\prec\circ_j\prec'$, that is it is either inverse for both $\prec\circ_j\prec'$ and shifted $\prec'$ or not.
\hfill $\Box$

\ 
 The same arguments  give the commutativity property.
 
  \begin{lemma}\label{comm} For a positive integer $d$ and nonnegative integers $k_0, k_1$,  
  
  
  there holds
 \[ (\prec\circ_{k_0}\prec')\circ_{k_0+m+k_1}\prec''= (\prec\circ_{k_0+d+k_1}\prec'')\circ_{k_0}\prec'.
 \]
 \end{lemma}
 
 $\Box$
 
 \

 \section{Operads of higher Bruhat orders}\label{op-higher-bruhat}

\subsection{Two operads} Fix a positive integer $d$. In this Section we will define two planar operads. The first one, 
to be denoted $HB_0^d$ and called {\em the small ($d$-th) Bruhat operad}, and the second one, to be denoted $HB^d$ and called {\em the big ($d$-th) Bruhat operad}. 

The operad $HB_0^d$ will be a suboperad of $HB^d$.

\subsection{Small Bruhat operad.}\label{small-br-op} Denote $HB_0^d(n) := B(nd,d)$. 


Define  composition maps
 
 \begin{equation}\label{compoper0}
 \gamma: HB_0^d(n)\times HB_0^d(m_1)\times \cdots \times HB_0^d(m_n)\to  HB_0^d(m_1+\cdots +m_n)
 \end{equation}
 by the rule
 
 \
 
 \[
\gamma(b_0; b_1, \ldots, b_n) = 
\]
 \[
 = (((b_0\circ_{0}b_1)\circ_{dm_1}\cdots )\circ_{dm_1\cdots + dm_{n-1}}b_n).
 \]

 \
 
 \begin{theorem}\label{main-thm}
The collection $HB_0^d(n)$,  $n\ge 1$, with the compositions  {\em (\ref{compoper0})}  
forms a planar operad $HB_0^d$. 
\end{theorem}

{\em Proof}.  Follows from Proposition \ref{inv-prop},   
Lemmata \ref{assoc},  \ref{comm}, and \ref{compos-ins}.   
 $\Box$
 
\subsection{Example: $d = 1$.} We have $HB_0^1(n) = B(n,1) = S(n)$; the insertion maps 
$$
\circ_j:\ S(m)\times S(n) \lra S(m+n-1)
$$
are described in \ref{ins-perm}. The resulting operad $HB_0^1$ appears in \cite{GKT}, 
Example 2.5, p. 21. 

\

\subsection{Big Bruhat operad.}\label{big-br-op} 

Denote by $HB^d(n)$ an infinite set whose elements are triples

\[
\alpha = (m, b, \mathbf{k})
\]
where $m$ is an integer $m > d$, $b\in B(m,d)$, and $\mathbf{k}$ is a type with 
$|\mathbf{k}| = m$ such that the molecule $M(\mathbf{k})$ has $n$ nuclei, each of them containing $d$ protons.

\

For notations see \ref{operad-molecules}.

\  
 
 Define  composition maps
 
 \begin{equation}\label{compoper}
 \gamma: HB^d(n)\times HB^d(m_1)\times \cdots \times HB^d(m_n)\to  HB^d(m_1+\cdots + m_n)
 \end{equation}
 by the following rule. Let $\alpha_0 = (m_0, b_0, \mathbf{k}_0)\in HB^d(n)$, 
 $\alpha_i = (m'_i, b_i, \mathbf{k}_i)\in HB^d(m_i), 1\leq i\leq n$. Let
 \[
 \mathbf{k}_0 = (k_0,d,k_1,\ldots,d,k_n)
 \] 
Then
\[
\gamma(\alpha_0; \alpha_1, \ldots, \alpha_n) = 
\]
 \[
 = (((b_0\circ_{k_0}b_1)\cdots )\circ_{k_0+m_1+k_1+\cdots +k_{n-2}+m_{n-1}}b_n).
 \]
 
 \
 

\begin{theorem}\label{main-thm-b}
The collection $HB^d(n)$ $n\ge 1$ with the compositions  {\em (\ref{compoper})}  forms a planar operad $HB^d$. We have an embedding of operads  
$HB_0^d\subset HB^d$. 
\end{theorem}
{\em Proof}.  Follows from Proposition \ref{inv-prop},   
Lemmata \ref{assoc},  \ref{comm}, and \ref{compos-ins}. 
$\Box$

\

For $d\le 2$, due to Corollary \ref{compat}
these planar operads are compatible with the $d$-th Bruhat order.

\

\

\centerline{Part II. MULTIPLICATION}

\section{Operads with multiplication, cosimplicial sets, shifted Poisson}\label{mult-cos-shifted}

\subsection{} 

Let $\mathcal{O} = \{\mathcal{O}(n), n\geq 0\}$ be a planar operad in the category $Sets$ of 
sets\footnote{or more generally in a tensor category whose objects are some structured sets, like topological spaces or abelian groups.}. Recall that this means that we are given a collection of 
sets $\mathcal{O}(n)$ together with composition maps
\[
\gamma: \mathcal{O}(n)\times \mathcal{O}(m_1)\times\ldots\times \mathcal{O}(m_n)\to 
\mathcal{O}(\sum_i m_i)
\]
and a unit element $1\in\mathcal{O}(1)$ satisfying some identities. 

\


Recall that defining an operadic composition is equivalent to defining a family of {\em insertion, or blowing up, maps}
\begin{equation}\label{ins}
\circ_j:\ \mathcal{O}(n)\times\mathcal{O}(m)\to \mathcal{O}(n+m-1),\ 0\leq j\leq n - 1,
\end{equation}
satisfying some simple identities, cf. \ref{part-comp}. 

\subsection{Operads with multiplication} Following \cite{MS}, Def. 3.1, let us call an {\em operad with multiplication} a 
planar operad $\mathcal{O}$ equipped with  elements 
$e\in \mathcal{O}(0), \mu\in \mathcal{O}(2)$ such that  
\[
\mu\circ_0\mu = \mu\circ_1\mu 
\]
and 
\[
\mu\circ_0 e = \mu\circ_1 e = 1\in \mathcal{O}(1).
\]

{\bf Change of notation.} 
We warn the reader that to be coherent with our previous notation the numbering of insertions will be {\it shifted by $1$} from 
the numbering of \cite{MS}:  they will start from $\circ_0$, and 
not from $\circ_1$ as in \cite{MS}.

\

Below we will show that our small and big Bruhat operads are in fact operads with multiplication. 

\subsection{Relation to cosimplicial sets} It is proven in \cite{MS} that an operad with multiplication $\mathcal{O} = \{\mathcal{O}(n)\}$ 
gives rise to a cosimplicial set $X  = \{ X^n\}$ whose set of cosimplices $X^n$ is 
$\mathcal{O}(n)$, see \ref{cos-sets} below.

Moreover this cosimplicial set acquires a family of maps
\[
\cdot : X^m\times X^n \to X^{m+n}
\]
compatible, in an appropriate sense, with cofaces and codegeneracies, cf. (\ref{comp-cos-a}) - (\ref{comp-cos-c}) below. Finally it is equipped with a family of elements $e_n\in X^n$ 
such that the collection of singletons $Y = \{ e_n\}$ is a cosimplicial subset of $X$. 

We call such structures {\em cosimplicial sets with multiplication}, and in such a way we get an equivalence of the category of operads with 
multiplication and the category of cosimplicial sets with multiplication, cf.  Thm. \ref{equiv-ops-cos}  
below. 

\subsection{Operads with multiplication and shifted Poisson}

Families of objects with insertions (\ref{ins}) appeared first (in an additive situation) in \cite{G} under the name of {\em pre-Lie systems}, cf. {\em op. cit.}, Section 5. It is proven 
in {\em op. cit.}, Section 6 that a pre-Lie system gives rise to a graded Lie algebra, 
see \ref{graded-lie} below.
 

\


\subsection{The Bruhat operads admit a multiplication} Sections 3 and 4 contain the main results of this paper. Recall that we have introduced two kinds of operads: 
{\em small Bruhat operads} $\Bscr_d$ and {\em big Bruhat operads} $\BBscr_d$; the big one 
contains the small one as a suboperad.   
 
 In Section 3 we shew that $\Bscr_d$ are operads with multiplication see Thm. \ref{small-bruhat-mult}, 
 whereas in Section 4 we shew that $\BBscr_d$ are operads with multiplication, 
 see Thm. \ref{big-bruhat-mult}.

\section{General remarks on operads with multiplication}

\subsection{Operads with multiplication and $Ass$.} Let $Ass$ denote a planar operad with $Ass(n)$ being a singleton 
$\{ e_n\}$ for all $n$, with a unique collection of compositions $\gamma$ and 
$1 = e_1\in Ass(1)$ satisfying the operadic identities. 

\

An  operad with multiplication is the same as a planar operad $\mathcal{O}$ equipped with a morphism of planar operads 
\[
\nu:\ Ass\to \mathcal{O},
\]
cf. \cite{MS}, Rem. 3.2 (i).

 

\

Namely, given $\nu$ as above, we define $e:= \nu(e_0), \mu = \nu(e_2)$.

\

Abusing the notation we will denote by the same symbol $e_n$ the element 
$\nu(e_n)\in\CO(n)$. 

\subsection{Cosimplicial sets.}\label{cos-sets} Let $(\mathcal{O}, e, \mu)$ be an operad with multiplication. Following \cite{MS}, Section 3, (10), p. 13, we assign to it a cosimplicial set 
$C(\CO) = C(\mathcal{O},\mu,e) = \{ C^n(\CO)\}$ as follows. 

\

We set $C^n(\CO) := \mathcal{O}(n)$. Define the cofaces $d^i: C^n(\CO)\to C^{n+1}(\CO)$ as follows:

\

$d^0(x) = \mu\circ_1 x; d^i(x) = x\circ_{i-1}\mu, 1 \leq i \leq n$, 
$d^{n+1}(x) = \mu\circ_0 x$

\



Define the codegeneracies $s^i: C^n(\CO)\to C^{n-1}(\CO)$ by $s^i(x) = x\circ_{i}e$.

\subsection{Multiplications $x\cdot y$.} Let $(\mathcal{O}, e, \mu)$ be an operad with multiplication. 
We have operadic compositions 
\[
\gamma_{2;m,n}:\ \mathcal{O}(2)\times\mathcal{O}(m)\times \mathcal{O}(n)\to 
\mathcal{O}(m+n), 
\]
whence the maps
\[
\cdot :\ \mathcal{O}(m)\times \mathcal{O}(n)\to \mathcal{O}(m+n),\
\]
\[ 
x\cdot y = \gamma_{2;m,n}(\mu;x,y) 
\]

\subsection{Claim.} {\em The multiplications $\cdot$ are associative.}

\

Cf. \cite{GV}, Prop. 2, (5).

\

They induce associative multiplications on the corresponding cosimplicial set 

$X = F(\mathcal{O},\mu,e)$, i.e. a family of maps

\begin{equation}\label{mult-cos}
\cdot:\ X^m\times X^n \to X^{m+n}
\end{equation} 

which satisfy the following compatibilities with cofaces and codegeneracies:

\begin{equation}\label{comp-cos-a}
d^i(x\cdot y) = (d^ix)\cdot y\ \text{if}\ i\leq m,\ = x\cdot d^{i-m}y\ \text{if}\ i>m;
\end{equation}

\begin{equation}\label{comp-cos-b}
s^i(x\cdot y) = (s^ix)\cdot y\ \text{if}\ i\leq m-1,\ = x\cdot s^{i-m}y\ \text{if}\ i\geq m
\end{equation}

("Leibnitz rules"), and

\begin{equation}\label{comp-cos-c}
(d^{m+1}x)\cdot y = x\cdot d^0y.
\end{equation}

cf. \cite{MS}, Def. 2.1 and Rem. 3.2 (ii).

\

So we have two kinds of multiplications:  $\cdot$ and $\circ_j$. \

\

Moreover, in each $X^n$ we have a distinguished element $e_n\in X_n$ such that
\begin{equation}\label{units-a}
d^i(e_n) = e_{n+1},\ s^i(e_n) = e_{n-1},
\end{equation}
and
\begin{equation}\label{units-b}
e_n\cdot e_m = e_{n+m}.
\end{equation}

\subsection{Claim} {\it For all $m, n, i$ 
\begin{equation}
e_m\circ_i e_n = e_{m+n-1}
\end{equation}}

\subsection{Definition.} {\em A cosimplicial set with multiplication} is a cosimplicial 
set $X$ equipped with a family of multiplications (\ref{mult-cos}) satisfying the identities (\ref{comp-cos-a}) - (\ref{comp-cos-c}) and with a family of elements $e_n\in X^n$  
satisfying the identities (\ref{units-a}) and (\ref{units-b}). 

\

Thus we have defined a functor 
\[
F: Opmult \to Cosmult
\] 
from the category of planar operads with 
multiplication to the category of cosimplicial sets with multiplication.

\subsection{Theorem.}\label{equiv-ops-cos} {\em The functor $F$ is an equivalence of categories.} $\Box$ 

\subsection{Additive setup: operation $\circ$.}
From now on till the end of this Section we will deal with a planar operad  
$\mathcal{V} = \{V(n)\}$  
in an additive tensor category. Following \cite{G}, Section 6 define operations 
\[
\circ:\ V(m)\otimes V(n) \to V(m+n-1)
\]
by 
\[
x_m\circ x_n = \sum_{i=0}^{m-1} (-1)^{(n+1)i} x_m\circ_i x_n.
\]
\
For example: 

\subsection{Claim}
\[
e_m\circ e_n=
\begin{cases}
me_{m+n-1} & \text{if} \, n \, \text{is} \, \text{odd},\cr
e_{m+n-1} & \text{if} \, n \,  \text{is} \, \text{even}, \, m\, \text{is} \, \text{odd},\cr
0 & \text{if} \, m\, \text{and} \, n\, \text{are}\, \text{even}.
\end{cases}
\]

\subsection{A graded Lie algebra}\label{graded-lie}
Let us introduce a bracket 
\[
[ , ]: V(m)\otimes V(n) \to V(m+n-1)
\]
by
\[
[x_m, x_n] = x_m\circ x_n - (-1)^{(m-1)(n-1)} x_n\circ x_m.
\]
Let us consider the graded space $V[1] = \oplus_n V(n)$, with $V(n)$ sitting in degree $n+1$.

The results of \cite{G}, Section 6 say that $V[1]$ with the above bracket $[ , ]$  
becomes a graded Lie algebra.

\

{\em Warning.} The operation $\circ$ is not associative, so the statement about a Lie algebra is not automatic.

\subsection{Switching multiplications on.} The structures below appear first in \cite{G}, Section 7 for the Hochschild complex of an associative ring, and  
for an arbitrary operad with multiplication (in an additive category)  in \cite{GV}, 
Section 1.  

\ 

Suppose in addition to the previous assumptions that $\mathcal{V}$ is an operad with multiplication. Then we get two news:

\

(a) The collection $V^\bullet = \{ V(n), n\geq 0\}$ becomes a cosimplicial set, whence we have  differentials
\[ 
d: V(n) \to V(n+1)
\] 
defined as usual as alternating sums of cofaces, 
\[
d = \sum_{i=0}^{n+1} (-1)^id^i.
\]
Similarly we may define 
\[ 
s: V(n) \to V(n-1)
\]
by
\[
s = \sum_{i=0}^{n-1} (-1)^is^i.
\]

\

(b) We have associative multiplications $\cdot: V(n)\otimes V(m) \to V(n+m)$.

\

{\em Warning.} The multiplication $\cdot$ is not commutative in general but it is commutative up to a homotopy, see Theorem \ref{homot-g} below. 

\

In particular we have an element $\pi = 1\cdot 1\in V(2)$. 

\subsection{Proposition.} (i) $\pi = d(1) = \mu$. (ii) $V = \oplus_n V(n)$ {\em equipped with the differential $d$ is a DG  associative algebra.}

\

See \cite{GV}, Section 1, Prop. 2 (5).

\subsection{Proposition} {\em We have}

(a) {\em $d(e_n) = 0$ if $n$ is even, and $d(e_n) = e_{n+1}$ if $n$ is odd.}

(b) {\em $s(e_n) = 0$ if $n$ is even, and $s(e_n) = e_{n-1}$ if $n$ is odd.}

\

Let us denote by $H^\cdot(V)$ the cohomology of $V = \oplus_n V(n)$ with respect to $d$.

\subsection{Proposition.} {\em We have
\[
x_m\cdot x_n = (\pi\circ_0 x_m)\circ_{m}x_n.
\]
for all $m, n$.} 

\

Cf. \cite{G}, Section 7, (22).

\

Now we have a fundamental

\subsection{Theorem.}\label{homot-g} {\em For all $x_m\in V(m), x_n\in V(n)$  
\begin{equation}\label{homot-g-relation}
- d(x_m\circ x_n) + (-1)^{n-1}dx_m\circ x_n + x_m\circ dx_n = 
(-1)^{n-1}(x_n\cdot x_m - (-1)^{mn}x_m\cdot x_n)
\end{equation}}

Cf. \cite{G}, Section 7, Thm. 3; \cite{GV}, Section 1, (9). 

\subsection{Homotopy Poisson} Moreover a trilinear operation
\[
h: V\otimes V\otimes V\to V[-2]
\]
is introduced in  \cite{GV} such that for all $x_m\in V(m), x_n\in V(n), x_k\in V(k)$
\[
[x_m, x_n\cdot x_k] - [x_m, x_n]\cdot x_k - (-1)^{m(n+1)}x_n\cdot [x_m, x_k] = 
\]
\[
= (-1)^{m+n}(dh(x_m\otimes x_n\otimes x_k) - hd(x_m\otimes x_n\otimes x_k)),
\]
see {\em op. cit.} (8).

\


\subsection{Corollary.} {\em The Lie bracket $[ , ]$ induces a Lie bracket on  
$H^\cdot(V)$.}

\    

It follows that $H^\cdot(V)$ is a $1$-shifted Poisson, or Gerstenhaber, algebra.

\

\subsection{Definition} Let us call {\em a Gerstenhaber - Voronov (GV) algebra} in an additive tensor category 
$\mathcal{A}$ a complex $V^\bullet$ in $\mathcal{A}$ with a differential 
$d$ of degree $1$ equipped with three operations: 

(i) a multiplication of degree $0$
\[
\cdot: V^\bullet\otimes V^\bullet \to V^\bullet
\]
making $V^\bullet$ an associative DG algebra; 

\

(ii) an operation of degree $-1$  
\[
\circ: V^\bullet\otimes V^\bullet \to V^\bullet[-1]
\]
such that if we define a bracket 
\[
[ , ]: V^\bullet\otimes V^\bullet \to V^\bullet[-1]
\]
by
\[
[x_n, x_m] = x_n\circ x_m - (-1)^{(n-1)(m-1)} x_m\circ x_n
\]
then $(V^\bullet[1], [ , ][1])$ becomes a graded Lie algebra;

\

(iii) an operation of degree $-2$
\[
h: V^\bullet\otimes V^\bullet\otimes V^\bullet \to V^\bullet[-2]
\]
such that
\[
[x_m, x_n\cdot x_k] - [x_m, x_n]\cdot x_k - (-1)^{m(n+1)}x_n\cdot [x_m, x_k] = 
\]
\[
= (-1)^{m+n+1}(dh(x_m\otimes x_n\otimes x_k) - hd(x_m\otimes x_n\otimes x_k)).
\]

\

(iv) The identity
\[
- d(x_m\circ x_n) + (-1)^{n-1}dx_m\circ x_n + x_m\circ dx_n = 
(-1)^{n-1}(x_n\cdot x_m - (-1)^{mn}x_m\cdot x_n)
\]
should hold. $\Box$

\

If $V^\bullet$ is a GV algebra then its cohomology $H^\bullet(V^\bullet)$ will be a Gerstenhaber algebra. 

\subsection{GV algebras vs BV algebras} The formula (\ref{homot-g-relation}) resembles 
but should not be confused with the defining relation  of a BV algebra.

\

Recall (cf. for example \cite{S}, Part II, Def. 2.1) that a {\em Batalin - Vilkovisky (BV) algebra} is a graded object $B^\bullet = \{ B^n\}$ equipped with a differential of degree $1$, a graded commutative and associative multiplication $\cdot$ of degree $0$, and a Lie bracket $[ , ]$ of degree $-1$ such that if we forget $d$, $B^\bullet$ is a Gerstenhaber algebra, and

\

for all $a\in B^m, b\in B^n$
\begin{equation}\label{bv-relation}
d(a\cdot b) - da\cdot b - (-1)^ma\cdot db  = (-1)^m [a, b], 
\end{equation}
compare this with (\ref{homot-g-relation}).

\


\section{Small Bruhat operads and multiplicative structures (Ursa Minor)}

\subsection{Insertions of higher Bruhat orders} Let $d$ be a positive integer. Let 
$B(m,d)$ denote the set of $d$-th Bruhat orders for the discrete Grassmanian 
$G(n,d)$ of subsets of cardinality $d$ in $[n] = \{1,\ldots,n\}$, cf. \cite{MaS}. 
Recall the insertion operations 
\begin{equation}\label{ins-bruhat}
\circ_j:\ B(m,d)\times B(n,d)\to B(m+n-d,d), 0\leq j\leq m+n-d,
\end{equation}
introduced in \ref{ins-higher-b}. 

\  

\subsection{The small operad of higher Bruhat orders}
Based on these operations,  a planar operad in $Sets$, to be denoted  
\newline $\mathcal{B}_d = \{\mathcal{B}_d(n)\}$ here,  has been introduced in 
\ref{small-br-op}, called {\em the small 
Bruhat operad}. By definition $\mathcal{B}_d(n) = B(nd,d)$.

The operadic compositions are defined in \ref{small-br-op}. Recall this formula: the map 
\begin{equation}
\gamma: \CB_d(n)\times \CB_d(m_1)\times \cdots \times \CB_d(m_n)\to  \CB_d(m_1+\cdots +m_n)
\end{equation}
is given by
\begin{equation}\label{gamma-formula}
\gamma(b_0; b_1, \ldots, b_n) 
 = (((b_0\circ_{0}b_1)\circ_{dm_1}\cdots )\circ_{dm_1\cdots + dm_{n-1}}b_n)
 \end{equation}
 for $b_0\in B(nd,d), b_i\in B(m_id, d), 1\leq i\leq n$.
 
 {\em Warning}: the insertion operations corresponding to compositions $\gamma$ should not 
 be confused with the operations (\ref{ins-bruhat}): we see from (\ref{gamma-formula}) that their numerotation is multiplied by $d$.  

\

Our aim in this Section 
will be to equip $\mathcal{B}_d$ with a structure of an operad with multiplication. 

\subsection{The case $d=1$.}  the case $d = 1$. The set $B(n,1)$ is the symmetric group $S(n) = Aut([n])$. It is equipped with the classical weak Bruhat order for which 
the minimal element is the identity permutation $e_n\in S(n)$, and the maximal one is the 
permutation of the maximal length $(n,n-1, \ldots, 2,1)$.

\

We set $[0] = \emptyset$, $\mathcal{B}_1(0) = S(0) = Aut([0])$; it is the singleton 
with the unique element $e_0$.

\

The insertion operations are
\[
\circ_j: S(n)\times S(m) \to S(n+m-1),
\]
see \ref{ins-perm}.

\

They make perfect sense for $m = 0$ as well, whence the maps 
\[
s^i := \circ_j: S(n) = S(n)\times S(0) \to S(n-1),\ 0\leq j \leq n-1,
\]
which are the codegeneracies in the corresponding cosimplicial set.

\subsection{Example} For $n = 3$ we have
\[
s^1(123) = s^2(123) = (12),\ s^3(123) = e.
\]
Here $(123)$ denotes as usually the cyclic permutation, taking $1$ to $2$, $2$ to $3$, 
and $3$ to $1$.


\subsection{Theorem.} {\em The triple $(\mathcal{B}_1, e_0, e_2)$ is an operad 
with multiplication in $Sets$.}

\



\subsection{Linearization.}
We can linearize this operad. For a set $I$ let $\mathbb{Z}I$ denote the free abelian group 
with base $I$. We have an obious embedding of sets $I\subset \mathbb{Z}I$. 

Let 
$\mathbb{Z}\mathcal{B}_1 = \{\mathbb{Z}\mathcal{B}_1(n)\}$ denote the planar operad 
in the tensor category $Ab$ of abelian groups, with $\mathbb{Z}\mathcal{B}_1(n) = 
\mathbb{Z}S(n)$, and operadic compositions 
induced by the compositions in $\mathcal{B}_1$. Then $(\mathbb{Z}\mathcal{B}_1, e_0, e_2)$ 
will be an operad with multiplication in $Ab$.

\

It gives rise to a GV algebra in $Ab$, to be denoted by $\BZ\mathbb{S}^\bullet$, with 
$\BZ\mathbb{S}^n = \mathbb{Z}S(n)$. 

\


\subsection{Cohomology} 


\begin{theorem}
We have \(H^0({\mathbb Z}S^\bullet)={\mathbb Z}e_0\), and
\(H^n({\mathbb Z}S^\bullet)=0\) for \(n>0\). 
\end{theorem}

\begin{proof}
Set \(C^n={\mathbb Z}S(n)\), and write
\(\delta^n=\sum_{i=0}^{n+1}(-1)^i d^i:C^n\to C^{n+1}\).  For \(n\geq 1\)
define the contracting homotopy \(h_n:C^n\to C^{n-1}\) on basis permutations
\(\sigma=(\sigma_1,\ldots,\sigma_n)\) by
\[
h_n(\sigma)=
\begin{cases}
(\sigma_2-1,\ldots,\sigma_n-1),& \sigma_1=1,\\
0,& \sigma_1\neq 1.
\end{cases}
\]
For \(n=1\), the first case is understood as \(e_0\in S(0)\).

We claim that, for every \(n\geq 1\),
\[
\delta^{n-1}h_n+h_{n+1}\delta^n=\operatorname{id}_{C^n}.
\]
It is enough to check this on a basis element
\(\sigma=(\sigma_1,\ldots,\sigma_n)\).  We use the following description of the
cofaces: \(d^0\sigma=(1,\sigma_1+1,\ldots,\sigma_n+1)\),
\(d^{n+1}\sigma=(\sigma_1,\ldots,\sigma_n,n+1)\), and, for \(1\leq i\leq n\),
\(d^i\sigma\) is obtained by splitting the \(i\)-th point of \(\sigma\).

First suppose that \(\sigma_1\neq 1\).  Then \(h_n\sigma=0\).  Among the cofaces
of \(\sigma\), only \(d^0\sigma\) begins with \(1\).  Hence
\(h_{n+1}d^0\sigma=\sigma\) and \(h_{n+1}d^i\sigma=0\) for \(1\leq i\leq n+1\).
Thus \((\delta^{n-1}h_n+h_{n+1}\delta^n)\sigma=\sigma\).

Now suppose that \(\sigma_1=1\), and put
\(\bar\sigma=h_n\sigma=(\sigma_2-1,\ldots,\sigma_n-1)\).  Then
\(h_{n+1}d^0\sigma=\sigma\), and, for \(1\leq i\leq n+1\),
\(h_{n+1}d^i\sigma=d^{i-1}\bar\sigma\).  Therefore
\[
h_{n+1}\delta^n\sigma
=\sigma+\sum_{i=1}^{n+1}(-1)^i d^{i-1}\bar\sigma
=\sigma-\sum_{i=0}^{n}(-1)^i d^i\bar\sigma
=\sigma-\delta^{n-1}\bar\sigma.
\]
Since \(\bar\sigma=h_n\sigma\), adding \(\delta^{n-1}h_n\sigma\) gives \(\sigma\).
This proves the claim.

If \(z\in C^n\), \(n>0\), and \(\delta^n z=0\), then the claim gives
\(z=\delta^{n-1}h_nz\).  Thus every positive-degree cocycle is a coboundary, so
\(H^n({\mathbb Z}S^\bullet)=0\) for \(n>0\).  In degree zero,
\(C^0={\mathbb Z}e_0\) and
\(\delta^0(e_0)=d^0(e_0)-d^1(e_0)=e_1-e_1=0\).  There is no incoming
differential into \(C^0\), hence \(H^0({\mathbb Z}S^\bullet)={\mathbb Z}e_0\).
\end{proof}

\subsection{Arbitrary $d$: the units.} Let now $d$ be arbitrary. Recall that the sets 
$B(m,d)$ are 
all ranged posets. Let $e_{m,d}\in B(m,d)$ denote the minimal element. It is caracterized by the property $Inv(e_{m,d}) = \emptyset$.

\subsection{Illustrations for $d = 2$} Elements of $B(m,2)$ are chains connecting 
the trivial permutation $e\in S(n) = B(n,1)$ with the maximal one (modulo some equivalence relation). They may be depicted by planar diagrams.

\

On  Fig. 3 below some unit elements $e_{m,2}\in B(m,2)$ are shewn.  

\begin{figure}[H]
\centering
\includegraphics[width=14cm, height=14cm]
{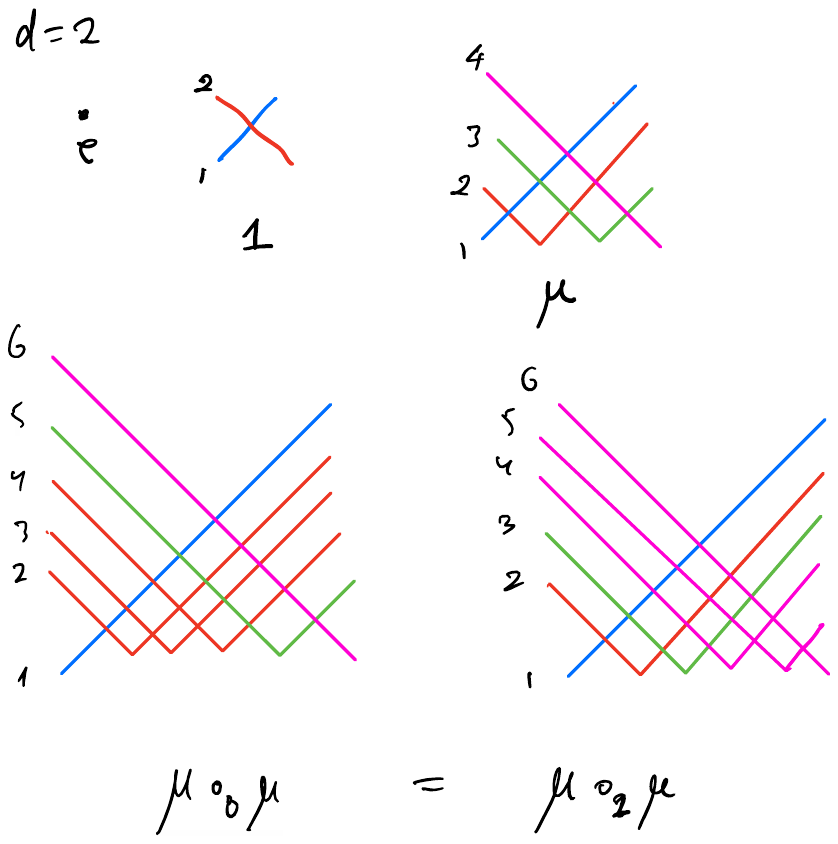}
\caption{Units in $B(m,2)$}
\end{figure}

On Figures 4, 5 an example of the products $x\cdot y$ and $y\cdot x$ is shewn; we see that 
the multiplication is not commutative.

\begin{figure}[H]
\centering
\includegraphics[width=14cm, height=14cm]
{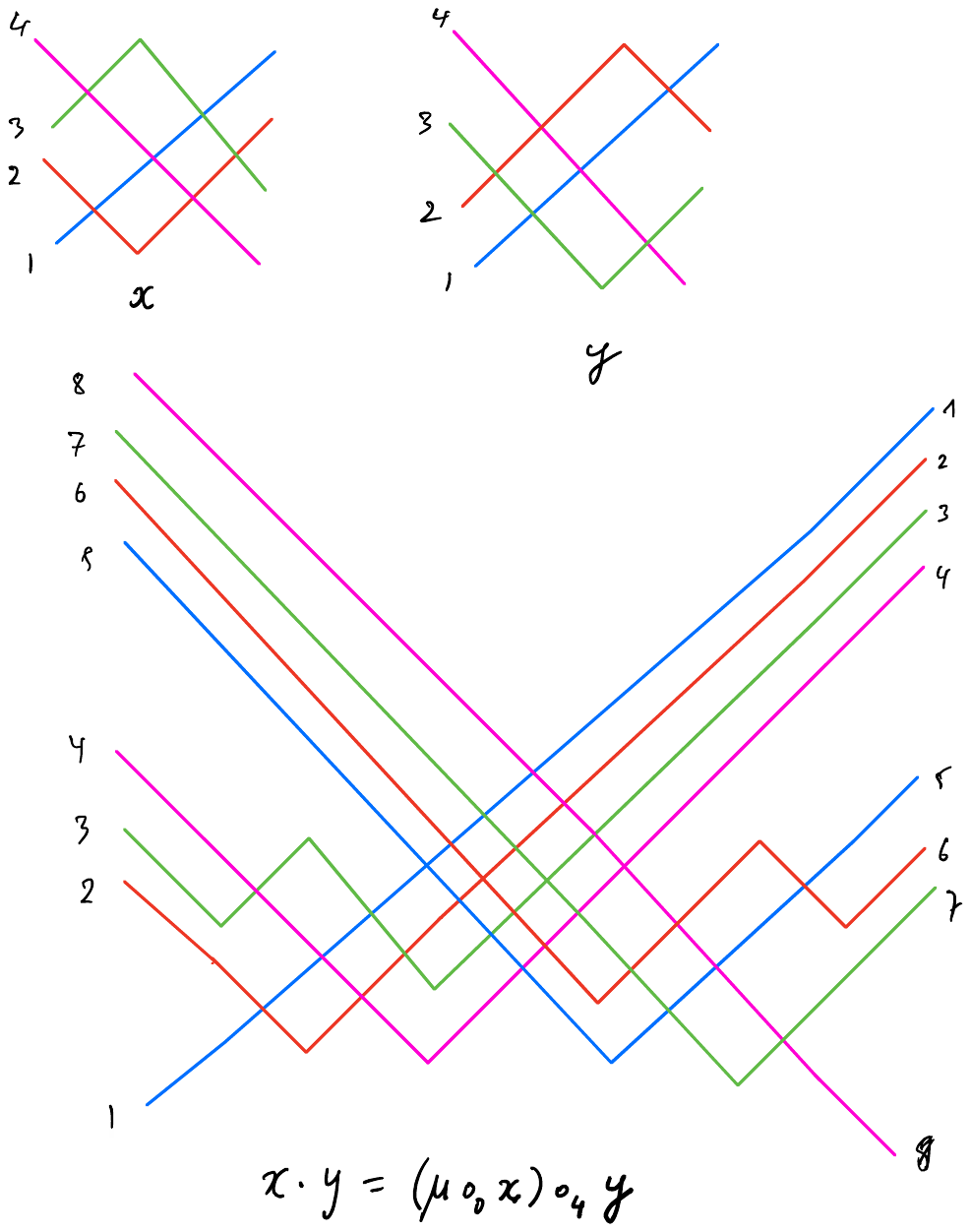}
\caption{Multiplication $x\cdot y$}
\end{figure}

\begin{figure}[H]
\centering
\includegraphics[width=14cm, height=14cm]
{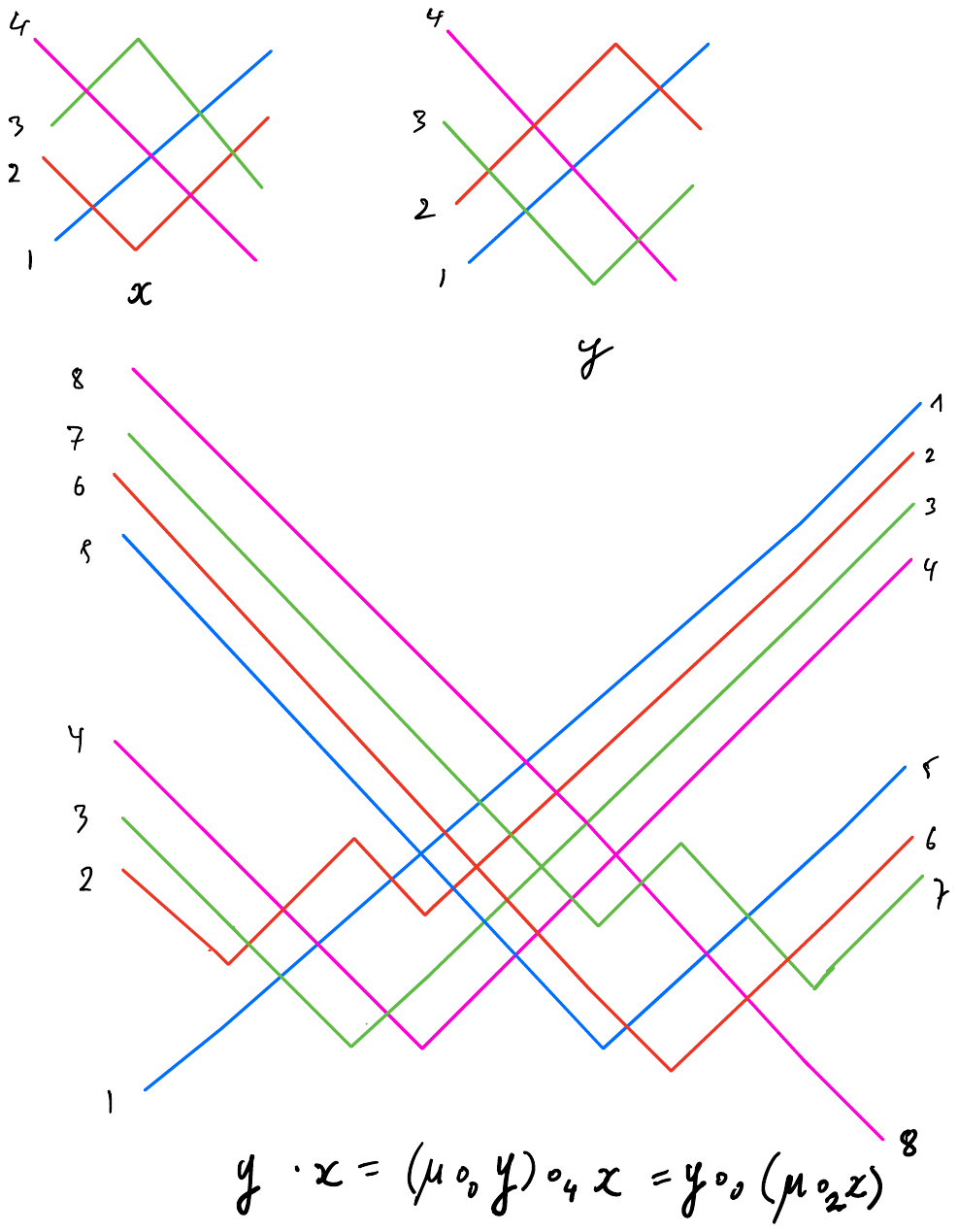}
\caption{Multiplication $y\cdot x$: different from $x\cdot y$}
\end{figure}

On Figure 6 all the cofaces of an element $x\in\CB_2(2)$ are shewn

\begin{figure}[H]
\centering
\includegraphics[width=14cm, height=14cm]
{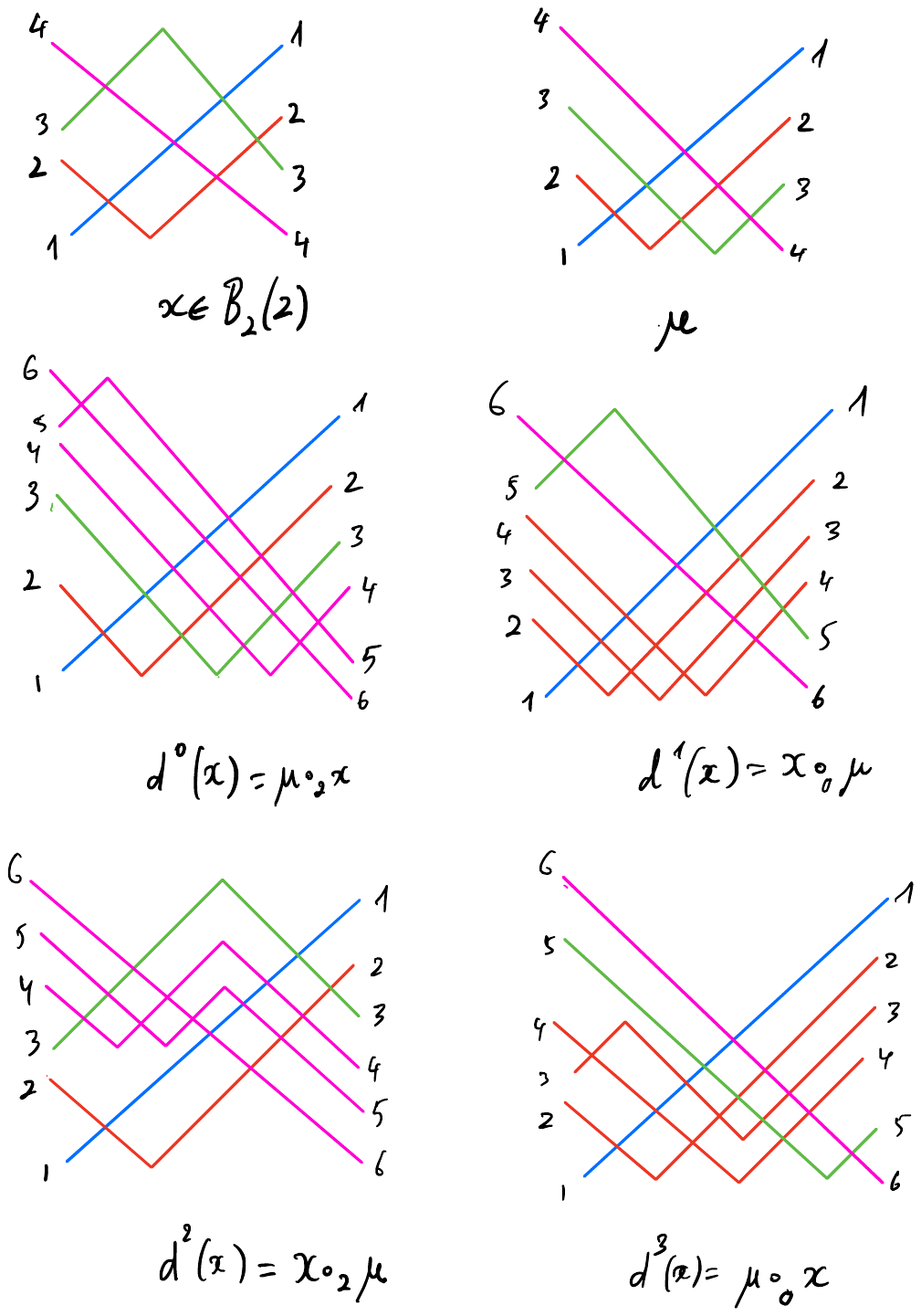}
\caption{Cofaces}
\end{figure}

On Figure 7 all the codegeneracies of the element $d^0(x)\in\CB_2(3)$ are shewn

\begin{figure}[H]
\centering
\includegraphics[width=14cm, height=14cm]
{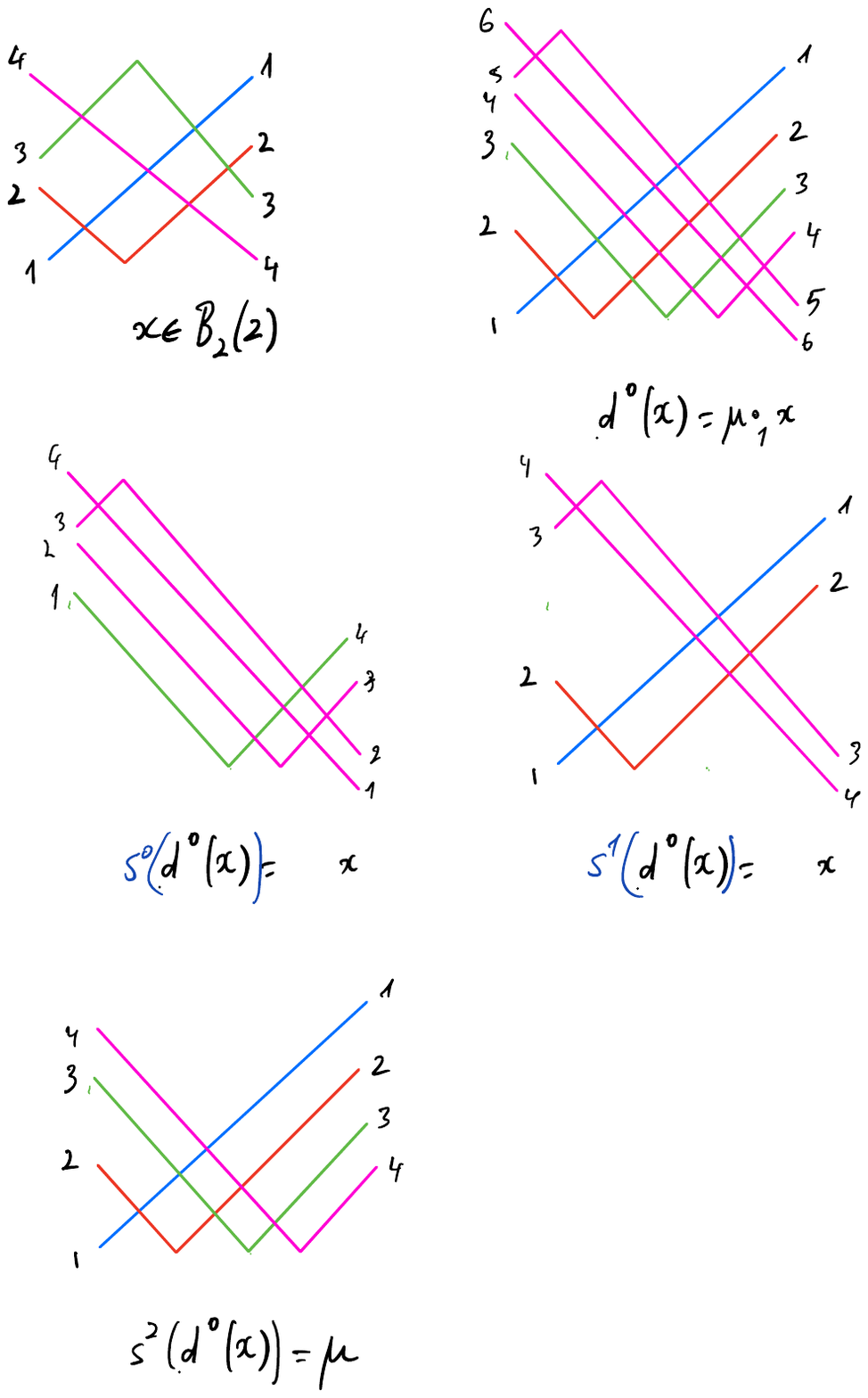}
\caption{Codegeneracies}
\end{figure}

\subsection{Theorem}\label{small-bruhat-mult} {\em Denote $e_n := e_{nd,d}\in B(nd,d) = \mathcal{B}_d(n)$. The triple $(\mathcal{B}_d, e_0, e_2)$ is an operad with multiplication.}

\

Linearizing as for the case $d = 1$ we get an operad $\mathbb{Z}\mathcal{B}_d$ in $Ab$.

\section{Big Bruhat operads and multiplicative structures (Ursa Major)}\label{ursa-major}

Similarly using the elements $e_{m,d}$ one defines a structure of an operad 
with multiplication on the {\em big} Bruhat operads from \ref{big-br-op}. 

\subsection{Master operad}
Let us recall the necessary definitions. Let us call {\em a type} a sequence 
\[
\mathbf{k} = (k_0,d,k_1,d,\ldots, k_{n-1},d,k_n)
\]
where $k_i$ are nonnegative integers. We denote $N = N(\mathbf{k}) := nd + \sum_{i=0}^n k_i$. 
 
Consider the interval  
\[
I := [N(\mathbf{k})] = \{1, \ldots, N(\mathbf{k})\};
\]
its elements will be called {\em particles}. 

Inside $I$ we have $n$ subintervals of length $d$ 
\[
I_j = k_j + jd + [d],\ 1\leq j\leq n
\]
called {\em nuclei},
whose elements are called {\em protons}. Elements of the complement 
$I\setminus\cup_{j=1}^n I_j$ are called {\em electrons}. 

This interval $I$ together with the above decomposition into $n$ nuclei of length $d$ and 
$N - nd$ electrons is called {\em the molecule of type $\mathbf{k}$} and denoted $M(\mathbf{k})$.

Let $\CM_d(n)$ denote the set of all molecules of all types containing 
$n$ nuclei of length $d$. It is shown in \ref{operad-molecules} that they form a planar operad 
$\CM_d = \{\CM_d(n), n\geq 0\}$, called the {\em Master operad}.  

We did not use the set $\CM_d(0)$ in {\em op. cit.} but our definition makes perfect sense 
for $n = 0$ as well and unlike the case of small Bruhat operad this set is not 
a singleton. Obviously all $\CM_d(n)$ are in bijection with $\BN^{n+1}$.

\subsection{The Master operad admits a multiplication} For any $n$ let $e_n\in \CM_d(n)$ denote the element with no electrons, i.e. 
of type 
\[
\mathbf{k}_0(n) := (0,d,0,\ldots, 0,d,0)
\] 
Then $(\CM_d, e_0, e_2)$ is an operad with multiplication.

\subsection{•} The $d$-th big Bruhat operad, to be denoted $\CBB_d = \{\CBB_d(n)\}$ here, 
is defined as follows, cf. \ref{big-br-op}. 

\

Elements of $\CBB_d(n)$ are couples $(b,\mathbf{k})$ where $\mathbf{k}\in \CM_d(n)$,  
$b\in B(m,d)$ with $m = N(\mathbf{k})$. 

\

Let $\CBB^0_d(n)\subset \CBB_d(n)$ denote the subset of couples of the form 
$(e_{m,d},\mathbf{k})$.

\subsection{Claim} {\em The collection $\{\CBB_d^0(n), n\geq 0\}$ is closed with respect to the operadic composition, so it forms a suboperad $\CBB^0_d\subset \CBB_d$. This suboperad 
is isomorphic to $\CM_d$.}

\

In particular we have the elements
\[
e_n = (nd, e_{nd,d}, \mathbf{k}_0(n))\in \CBB_d^0(n)\subset \CBB_d(n).
\]

\subsection{Theorem}\label{big-bruhat-mult} {\em The triple $(\CBB_d, e_0, e_2)$ is an operad with multiplication 
in $Sets$.}

\

\centerline{Part III. APPENDIX by DARIA POLIAKOVA}

\section{}
The purpose of this Appendix is to prove the a strengthening of Corollary~5.4.

Elements of \(B(n,d)\) carry two natural orders.  The first one is the
dominance order, defined by inclusion of inversion sets.  The second one is
the flip order, whose covering relations are the increasing flips, i.e. the
addition of one inversion.  Corollaries~5.3 and~5.4 show that insertion
respects the dominance order in the second and in the first argument,
respectively.  For Corollary~5.3 the corresponding statement for the flip
order is immediate: one inversion in the second argument becomes one shifted
inversion after insertion. For Corollary~5.4, an additional argument is required.

\begin{theorem}\label{thm:insertion-flip-order}
The insertion operations are monotone for the flip order.  Let
\(\prec, \widetilde{\prec}\in B(n,d)\), let
\(\prec',\prec''\in B(m,d)\), and let \(j\) be an admissible insertion place.  Then
\[
  \prec' \leq_{\mathrm{fl}} \prec''
  \quad\Longrightarrow\quad
  \prec\circ_j\prec' \leq_{\mathrm{fl}} \prec\circ_j\prec'',
\]
and
\[
  \prec \leq_{\mathrm{fl}} \widetilde{\prec}
  \quad\Longrightarrow\quad
  \prec\circ_j\prec' \leq_{\mathrm{fl}} \widetilde{\prec}\circ_j\prec' .
\]
\end{theorem}

For the proof, recall the insertion construction from 5.8.  Let \(\prec\in B(n,d)\),
\(\prec'\in B(m,d)\), and \(0\leq j\leq n-d\).  The insertion
\(\prec\circ_j\prec'\in B(n+m-d,d)\) is defined by its inversion set.  Identify
\([n+m-d]=A\sqcup B\sqcup C\), where \(A=[j]\), \(B=j+[m]\), and
\(C=j+m+[\,n-(j+d)\,]\), and let \(D=j+[d]\subset B\).  We identify \([m]\)
with \(B\) by translation by \(j\), and \([n]\) with \(A\sqcup D\sqcup C\) by the
monotone bijection \(\phi:[n]\to A\sqcup D\sqcup C\).  For
\(L\in\binom{A\sqcup B\sqcup C}{d+1}\), consider three cases.  If \(L\subset B\),
then \(L\in\operatorname{Inv}(\prec\circ_j\prec')\) if and only if
\(L-j\in\operatorname{Inv}(\prec')\).  If \(L\subset A\sqcup D\sqcup C\), then
\(L\in\operatorname{Inv}(\prec\circ_j\prec')\) if and only if
\(\phi^{-1}(L)\in\operatorname{Inv}(\prec)\).  Finally, in the remaining case,
define \(\bar L\) by moving the elements of \(L\) that are in \(B\setminus D\) into
the rightmost available positions of \(D\).  Then
\(L\in\operatorname{Inv}(\prec\circ_j\prec')\) if and only if
\(\phi^{-1}(\bar L)\in\operatorname{Inv}(\prec)\).

\begin{proof}
If \(\prec''\) is obtained from \(\prec'\) by an increasing flip adding
\(F\in\binom{[m]}{d+1}\), then, by the case \(L\subset B\) in the definition of
insertion, \(\prec\circ_j\prec''\) is obtained from \(\prec\circ_j\prec'\) by the
single increasing flip adding \(j+F\subset B\); hence the monotonicity in the second
argument follows by iteration.  It remains to prove the monotonicity in the first
argument, and again it is enough to assume that \(\widetilde{\prec}\) is obtained
from \(\prec\) by one increasing flip, that is,
\(\operatorname{Inv}(\widetilde{\prec})=\operatorname{Inv}(\prec)\sqcup\{K\}\),
\(K\in\binom{[n]}{d+1}\). We then need to show that there is a sequence of increasing
flips from \(\prec\circ_j\prec'\) to \(\widetilde{\prec}\circ_j\prec'\).

The strategy of our proof is as follows.
\begin{enumerate}
    \item We first describe the set of flips that we need, and identify them with
    \(q\)-subsets of some set \(H\).
    \item We then construct an element of \(B(|H|,q)\) by its inversion set, and check
    the Z-criterion.
    \item We then choose a lift to an admissible order in \(A(|H|,q)\), and show that
    the flips can be performed in this order, i.e. that each newly obtained set satisfies
    the Z-criterion.
\end{enumerate}

\paragraph{Step 1.}
We need to describe the inversions that are in
\(\operatorname{Inv}(\widetilde{\prec}\circ_j\prec')\) but not in
\(\operatorname{Inv}(\prec\circ_j\prec')\); these are our flips of interest. Put
\(K^\phi:=\phi(K)\subset A\sqcup D\sqcup C\), and set formally
\(\overline{K^\phi}=K^\phi\). The flips of interest correspond to
\[
L \in \binom{A \sqcup B \sqcup C}{d+1}
\]
such that \(\bar L = K^\phi\). Let \(Q\) be the terminal segment of \(D\) occupied by
elements of \(K^\phi\). If the last element of \(D\) is not in \(K^\phi\), then \(Q\) is
empty and the only \(L\) with \(\bar L=K^\phi\) is \(K^\phi\) itself, so the statement
follows. We thus assume that \(Q\) is not empty.
To have \(\bar L=K^\phi\), we must have the elements of \(Q\) either already in \(Q\)
in \(L\), or have them in \(B\setminus D\), so that we move them to \(Q\) when forming
\(\bar L\). Set
\[
q = |Q|,
\qquad
H = Q \sqcup (B \setminus D).
\]
Then the flips of interest are in bijection with the \(q\)-subsets of \(H\). Indeed, let
\(K^{-}\) denote the elements of \(K^\phi\) before \(H\), and let \(K^{+}\) denote the
elements of \(K^\phi\) after \(H\). For any subset \(F \subset H\), set its extension
\[
\widetilde F := K^{-} \sqcup F \sqcup K^{+}.
\]
Then for every \(q\)-subset \(F\) of \(H\), \(\widetilde F\) is a \((d+1)\)-subset of
\(A \sqcup B \sqcup C\) which bar-projects onto \(K^\phi\).

\paragraph{Step 2.}
We now want to construct an element of \(B(|H|,q)\) by describing its inversion set.
Let \(U = \{u_1 < \cdots < u_{q+1}\}\) be a \((q+1)\)-subset of \(H\). The
\(q\)-subsets
\[
U_1 = U \setminus \{u_{q+1}\}, \ldots, U_{q+1} = U \setminus \{u_1\}
\]
are its packet, and their extensions
\[
\widetilde{U_i} = K^{-} \sqcup U_i \sqcup K^{+}
\]
are a contiguous segment in the packet of
\[
\widetilde U = K^{-} \sqcup U \sqcup K^{+}.
\]
By the previous step, all \(\widetilde{U_i}\) are in
\(\operatorname{Inv}(\widetilde{\prec}\circ_j\prec')\), but not in
\(\operatorname{Inv}(\prec\circ_j\prec')\). Since both
\(\operatorname{Inv}(\widetilde{\prec}\circ_j\prec')\) and
\(\operatorname{Inv}(\prec\circ_j\prec')\) satisfy the Z-criterion, it follows that
\(\operatorname{Inv}(\prec\circ_j\prec')\) contains either only the prefix of
\(P(\widetilde U)\) before the segment of the \(\widetilde{U_i}\)'s, or only the suffix
after that segment. Set \(\Theta\) to consist of all \(U\) such that, in the above
dichotomy, it is the suffix.

We now need to show that \(\Theta\) satisfies the Z-criterion. First suppose that
\(K^{-}\) is nonempty, and choose some \(a \in K^{-}\). A \((q+1)\)-subset \(U\)
belongs to \(\Theta\) if and only if
\[
\widetilde U \setminus \{a\}
=
(K^{-} \setminus \{a\}) \sqcup U \sqcup K^{+}
\]
belongs to \(\operatorname{Inv}(\prec\circ_j\prec')\), since
\((K^{-} \setminus \{a\}) \sqcup U \sqcup K^{+}\) goes after the segment of the
\(\widetilde{U_i}\)'s in \(P(\widetilde U)\). Now let
\(V = \{v_1 < \cdots < v_{q+2}\}\) be some \((q+2)\)-subset of \(H\). By the
reasoning above, \(V \setminus \{v_i\}\) belongs to \(\Theta\) if and only if
\[
\widetilde V \setminus \{a,v_i\} \in \operatorname{Inv}(\prec\circ_j\prec').
\]
In the packet of \(\widetilde V \setminus \{a\}\), the sets
\(\widetilde V \setminus \{a,v_i\}\) form a contiguous segment, so by the Z-criterion
for \(\operatorname{Inv}(\prec\circ_j\prec')\), either a beginning or a terminal segment
of them belongs to \(\operatorname{Inv}(\prec\circ_j\prec')\), proving the Z-criterion
for \(\Theta\). Finally, if \(K^{-}\) is empty, then \(K^{+}\) is nonempty, so one
chooses some \(b \in K^{+}\) and runs essentially the same argument with \(b\).

\paragraph{Step 3.}
We have shown that \(\Theta\) satisfies the Z-criterion; thus there exists an
admissible order on the \(q\)-subsets of \(H\) whose inversion set is \(\Theta\). Let
\[
W_1, W_2, \ldots, W_r,
\qquad
r = \binom{|H|}{q},
\]
be one such order. We claim that the flips \(\widetilde{W_i}\) can be performed in
that order, or, equivalently, that for any \(1 \leq t \leq r\), the set
\[
I_t = \operatorname{Inv}(\prec\circ_j\prec')
      \sqcup \{\widetilde{W_1}, \ldots, \widetilde{W_t}\}
\]
satisfies the Z-criterion. Let \(S\) be a \((d+2)\)-subset to be tested.

If \(P(S)\) does not contain any \(\widetilde{W_i}\), then the intersection of \(I_t\)
with \(P(S)\) is the same as the intersection of
\(\operatorname{Inv}(\prec\circ_j\prec')\) with \(P(S)\), which is either initial or
terminal. If \(P(S)\) contains exactly one \(\widetilde{W_i}\), then the intersection of
\(I_t\) with \(P(S)\) is either the same as for
\(\operatorname{Inv}(\prec\circ_j\prec')\) if \(i > t\), or the same as for
\(\operatorname{Inv}(\widetilde{\prec}\circ_j\prec')\) if \(i \leq t\), both of which are
either initial or terminal.

Finally, suppose that \(P(S)\) contains \(\widetilde{W_i}\) and \(\widetilde{W_\ell}\)
for \(i \neq \ell\). Then \(U = W_i \cup W_\ell\) is a \((q+1)\)-subset of \(H\). If
\(U \in \Theta\), then by definition the suffix of \(P(\widetilde U)\) after the
segment of the \(\widetilde{U_i}\)'s is in
\(\operatorname{Inv}(\prec\circ_j\prec')\), and the order
\(W_1,\ldots,W_r\) inverts the lexicographic order on the packet of \(U\). Thus the
already added elements in this packet form a terminal segment of this packet, and,
together with the suffix, we obtain the terminal segment of the whole packet \(P(S)\).
If \(U \notin \Theta\), then by definition the prefix of \(P(\widetilde U)\) before the
segment of the \(\widetilde{U_i}\)'s is in \(\operatorname{Inv}(\prec\circ_j\prec')\),
and the order \(W_1,\ldots,W_r\) restricts to the lexicographic order on the packet of
\(U\). Thus the already added elements in this packet form an initial segment of this
packet, and, together with the prefix, we obtain the initial segment of the whole packet
\(P(S)\).

This completes the proof.
\end{proof}


\

G.K.: Institute for Information Transmission Problems Russian Academy of Sciences,
Russian Federation; 
{koshevoyga@gmail.com}

V.S.: Institut de Math\'ematiques de Toulouse, 
Universit\'e Toulouse III Paul Sabatier, 119 Route de Narbonne, 31062 Toulouse, France; 
Kavli IPMU, UTIAS, 5-1-5 Kashiwanoha, Kashiwa, Chiba, 277-8583 Japan; 
{schechtman@math.ups-tlse.fr} 

D.P.: University of Hamburg, Bundesstrasse 55, 20146 Hamburg, Germany; 
{polydarya@gmail.com}

\end{document}